\documentclass[11pt,a4paper]{article}

\usepackage{graphicx}
\usepackage{amsmath}
\usepackage{amsfonts}
\usepackage{amssymb}
\usepackage{enumerate}
\usepackage{lscape}
\usepackage{longtable}
\usepackage{rotating}
\usepackage{color}
\usepackage{url}
\usepackage{subfigure}
\usepackage{rotating}
\newtheorem{theorem}{Theorem}[section]

\newtheorem{algorithm}{Algorithm}[section]

\newtheorem{conjecture}{Conjecture}[section]

\newtheorem{definition}{Definition}[section]

\newtheorem{lemma}{Lemma}[section]

\newtheorem{assumption}{Assumption}[section]

\newcommand{\beqn}[1]{\begin{equation}\label{#1}}
\newcommand{\eeqn}{\end{equation}}
\def\NN{\mathbb{N}}
\def\EE{\mathbb{E}}
\def\BBB{{\bf b}}

\newcommand{\calP}{\mathcal{P}}

\newcommand{\R}{\mathbb{R}}
\newcommand{\half}{\frac{1}{2}}

\DeclareMathOperator{\cond}{cond}

\newcommand{\Prob}{\operatorname{Prob}}

\newcommand{\proof}[1]{{\noindent\emph{Proof.}} #1 \hfill$\Box$ \linebreak}

\begin{document}
\title{Convergence of  trust-region methods \\ based on probabilistic models}
\author{
A. S. Bandeira\thanks{Program on Applied and Computational Mathematics,
Princeton University, Princeton, NJ 08544, USA ({\tt ajsb@math.princeton.edu}). Support for this author was provided by NSF Grant No. DMS-0914892.}
\and
K. Scheinberg\thanks{
Department of Industrial and Systems Engineering, Lehigh University,
Harold S. Mohler Laboratory, 200 West Packer Avenue, Bethlehem, PA 18015-1582, USA
({\tt katyas@lehigh.edu}). The work of this author is partially supported by NSF Grant DMS 10-16571, AFOSR Grant FA9550-11-1-0239, and  DARPA grant FA 9550-12-1-0406 negotiated by AFOSR.}
\and
L. N. Vicente%
\thanks{CMUC, Department of Mathematics, University of Coimbra,
3001-501 Coimbra, Portugal ({\tt lnv@mat.uc.pt}). Support for
this author was provided by FCT under grants PTDC/MAT/116736/2010 and PEst-C/MAT/UI0324/2011.}
}
\maketitle
\footnotesep=0.4cm
{\small
\begin{abstract}
In this paper we consider the use of probabilistic or random models within a classical
 trust-region  framework for optimization of deterministic smooth general nonlinear functions.
 Our method and setting differs from many  stochastic optimization approaches in two principal ways.
 Firstly, we assume that the value of the function itself can be computed without noise, in other words, that the function is deterministic.
 Secondly,  we use random models of higher quality than those produced by usual stochastic gradient methods. In particular,
 a first order model based on  random approximation of the gradient is required to provide sufficient quality of approximation
 with probability greater than or equal to  $1/2$. This is in contrast with stochastic gradient approaches, where the model is assumed to be ``correct''
only  in expectation.

As a result of this particular setting, we are able to prove convergence, with probability one, of a trust-region method which
is almost identical to the classical method. Hence we show that  a standard optimization framework can be used in cases when models are random and may or may not provide  good approximations, as long as ``good'' models are more likely than ``bad'' models. Our results are based on the use of properties of martingales. Our motivation comes from using random sample sets and interpolation models in derivative-free optimization. However, our framework is general and can be applied with any source of uncertainty in the model. We discuss various applications for our methods in the paper.
 \end{abstract}

\bigskip

\begin{center}
\end{center}
}


\section{Introduction}

\subsection{Motivation}

The focus of this paper is the analysis of a  numerical scheme that utilizes randomized models to minimize deterministic functions.  In particular, our motivation comes from algorithms for  minimization of so-called \emph{black-box} functions where values are computed, e.g., via simulations.  For such problems, function evaluations are costly and derivatives are typically unavailable and cannot be approximated.  Such is the setting of \emph{derivative-free optimization} (DFO), of which the list of applications---including molecular geometry optimization, circuit design, groundwater community problems, medical image registration, dynamic pricing, and aircraft design (see the references in~\cite{ARConn_KScheinberg_LNVicente_2009})
--- is diverse and growing. Nevertheless,
our framework is general and is not limited to the setting of derivative-free optimization.

There is a variety of evidence supporting the claim that randomized models can yield both practical and theoretical benefits for deterministic optimization.  A primary example is the recent success of stochastic gradient methods for solving large-scale machine learning problems.  As another example, the randomized coordinate descent method for large-scale convex deterministic optimization proposed in~\cite{YNesterov_2012} yields better complexity results than, e.g., cyclical coordinate descent.   Most contemporary randomized methods generate random directions along which all that may be required is some minor level of descent in the objective $f$.  The resulting methods may be very simple and enjoy low per-iteration complexity, but the practical performance of these approaches can be very poor.  On the other hand, it was noted in~\cite{RByrd_etal_2012} that the performance of stochastic gradient methods for large-scale machine learning improves substantially if the sample size is increased during the optimization process.
Within direct search, the use of random positive bases has also been recently
investigated~\cite{CAudet_JEDennis_2006,LNVicente_ALCustodio_2009} with gains
in performance and convergence theory for nonsmooth problems.
This suggests that for  a wide range of optimization problems, requiring a higher level of accuracy from a randomized model may lead to more efficient methods.  Thus, our primary goal is to design randomized numerical methods that do not rely on producing descent directions ``eventually'', but provide accurate enough approximations  so that in each iteration a sufficiently improving step is produced with high probability (in fact, probability greater than half is sufficient in our analysis).  We  incorporate these models into a trust-region framework so that the resulting algorithm is
able to work well in practice.

Our motivation originates with model-based  DFO methodology (e.g., see \cite{ARConn_KScheinberg_LNVicente_2009-paper,ARConn_KScheinberg_LNVicente_2009}) where local models of $f$ are built from function values sampled in the vicinity of a given iterate.  To date, most algorithms of this type have relied on sample sets that are generated by the algorithm steps or added in a deterministic manner.  A complex mechanism of sample set maintenance is necessary to ensure that the quality of the models is acceptable, while the expense of sampling the function values is not excessive.  Various approaches have been developed for this mechanism, which achieve different trade-offs for the number of sample points required, the computational expense of the mechanism itself, and the quality of the models.  One of the primary premises  of this paper is the assumption that using random sample sets can yield new and better trade-offs.  That is, randomized models can maintain a higher quality by using fewer sample points without complex maintenance of the sample set.  One example of such a situation is described in~\cite{ASBandeira_KScheinberg_LNVicente_2012}, where
 linear or quadratic polynomial models are constructed
from random sample sets. It is shown that one can
build such models, meeting a Taylor type accuracy with high probability, using significantly less sample
points than what is needed in the deterministic case, provided the function being
modeled has sparse  derivatives.

The framework considered by us in the current paper is sufficiently broad
to encompass any situation where the quality or accuracy of the trust-region
models is random. In particular, such models can be built directly using some form
of derivative information, as long as it is accurate with  certain probability.

\subsection{Trust-region framework}

The trust-region method introduced and analyzed in this paper is rather simple. At each iteration one solves
a trust-region subproblem, i.e., one minimizes the model within a trust-region ball.
Note that one does not know whether the model is accurate or not. If the trust-region
step yields a good decrease in the objective function relatively to the decrease in the
model and the trust-region radius is sufficiently small relatively to the size of the model gradient,
then the step is taken and the trust-region radius is possibly increased. Otherwise the step
is rejected and the trust-region radius is decreased.
We show that such a method  always drives the trust-region radius to zero.

Based on this property we show that, provided the  (first order) accuracy of the model occurs with probability
no smaller than $1/2$, possibly conditioned to the prior iteration history, then
the gradient of the objective function converges to zero with
probability one. Our proof technique relies on building random processes from the random
events defined by the models being or not accurate (conditioned to the past),
and then making use of their submartingale-like properties. We extend the theory
to the case when the models of sufficient  second order accuracy  occur with probability
no smaller than $1/2$. We show that a subsequence of the iterates drive a measure of
second order stationarity to zero with probability one. However, to demonstrate the $\lim$-type convergence to a second order stationary point
we need additional assumptions on the models.

\subsection{Notation}

Several constants are used in this paper to bound various quantities. These constants
are denoted by $\kappa$ with acronyms for the subscripts that are indicative of the quantities
that they are meant to bound.
We list their most used definitions here,
for convenience. The actual meaning of the constants will become clear when each
of them is introduced in the paper.

\vspace{1ex}
\begin{tabular}{ll}
$\kappa_{fcd}$ & ``fraction of Cauchy decrease'' \\
$\kappa_{fod}$ & ``fraction of optimal decrease''\\
$\kappa_{Lg}$ & ``the Lipschitz constant of the gradient of the function'\\
$\kappa_{Lh}$ & ``the Lipschitz constant of the Hessian of the function''\\
$\kappa_{L\tau}$ & ``the Lipschitz constant of the measure $\tau$ of second order stationarity of the function''\\
$\kappa_{ef}$ & ``error in the function value''\\
$\kappa_{eg}$ & ``error in the gradient''\\
$\kappa_{eh}$ & ``error in the Hessian''\\
$\kappa_{e\tau}$ & ``error in the $\tau$ measure''\\
$\kappa_{bhm}$ & ``bound on the Hessian of the models''\\
$\kappa_{bhf}$ & ``bound on the Hessian of the function''
\end{tabular}
\vspace{1ex}

This paper is organized as follows. In  Section~\ref{sec:dfo} we briefly describe existing methods for derivative-free optimization and provide an illustrative example
to motivate the use of random models. In  Section~\ref{sec:first-order}
we introduce the probabilistic models of the first order and the trust-region
method based on such models. The convergence of the method to first order
criticality points is proved in Section~\ref{sec:conv-1st-order}.
The second order case is addressed in Section~\ref{sec:2nd-order}. Finally,
in Section~\ref{sec:rand_models} we describe various useful random models that satisfy the conditions needed for convergence results in
Sections~\ref{sec:first-order} and~\ref{sec:2nd-order}.

\section{Methods of derivative-free optimization}\label{sec:dfo}

We consider in this paper the unconstrained optimization problem
\[
\min_{x \in \mathbb{R}^n} f(x),
\]
where the first (and second, in some cases) derivatives of the objective function $f(x)$ are
assumed to exist and be Lipschitz continuous.  However,
as it is considered in derivative-free optimization (DFO),
explicit evaluation of these derivatives is assumed to be impossible.
Derivative-free methods rely on sampling the objective function
at one or more points at each iteration. Some sample to explore directions,
others to build models.

\paragraph*{Directional methods.}

Direct-search methods of directional type were first developed using a single positive basis or a finite number of
them (see the surveys~\cite{TGKolda_RMLewis_VTorczon_2003}
and \cite[Chapter~8]{ARConn_KScheinberg_LNVicente_2009}).
The basic versions of these methods, like coordinate or compass search,
are inherently slow for problems of more than a few variables,
not only because they are not able to use curvature information and rarely reuse sample points,
but also because they rely on few directions.
They were shown to be globally convergent for smooth problems~\cite{VTorczon_1997}
and had their worst case complexity measured by global rates~\cite{LNVicente_2013}.

Not restricting direct search to a finite number of positive bases was soon
discovered to enhance practical performance.
Approaches allowing for an infinite number of positive bases were
proposed in~\cite{CAudet_JEDennis_2006,TGKolda_RMLewis_VTorczon_2003},
with results applicable to nonsmooth functions
when the generation is dense in the unit sphere
(see~\cite{CAudet_JEDennis_2006,LNVicente_ALCustodio_2009}).

On the other hand,
randomized stochastic methods recently became  a popular alternative to
direct-search methods. These methods also sample the objective function
along a certain direction, but instead of choosing a direction from
a positive basis, these methods select directions totally randomly. This often allows for
faster convergence because ``good'' directions are occasionally observed.
The random search approach introduced in \cite{JMatyas_1965} samples points from a Gaussian
distribution.  Convergence of an improved scheme was shown in~\cite{BTPoljak_1987}.
In \cite{YNesterov_2011}, Nesterov recently presented several derivative-free random search schemes and provided bounds for their global rates. Different improvements of these methods emerged in the latest literature, e.g.,~\cite{SGhadimi_GLan_2012}.
Although complexity results for both convex and nonsmooth nonconvex functions are available for randomized search,
the practical usefulness of these methods is limited by the fixed step sizes determined by the complexity analysis and,
 as in direct search, by the lack of curvature information.

\paragraph*{Model-based trust-region methods.}

Model-based DFO methods developed by
Powell~\cite{MJDPowell_1994,MJDPowell_2001,MJDPowell_2003,MJDPowell_2004},
and by Conn, Scheinberg, and Toint
\cite{ARConn_KScheinberg_PhLToint_1997a,ARConn_KScheinberg_PhLToint_1997b}, introduced a class of trust-region methods
that relied on  {\em  interpolation} or {\em regression} based quadratic approximations
of the objective function  instead of the usual Taylor series
quadratic   approximation. The regression-based method
was later successfully used in \cite{SCBillups_JLarson_PGraf_2013}
based on~\cite{ARConn_KScheinberg_LNVicente_2008b}.
In all cases the models are built based on  sample points  in
reasonable proximity to the current  best iterate.
The computational study of Mor\'e and Wild
\cite{JJMore_SMWild_2009} has shown that these methods are
typically significantly superior in practical performance to the other existing approaches
due to the use of  models  that effectively capture the  local curvature of the objective function.
While the model quality is undoubtedly essential for the performance of these methods,
guaranteeing sufficient  quality at all times is quite expensive computationally.
Randomized models, on the other hand, can offer a suitable alternative by providing a good quality
approximation with high probability.

\paragraph*{An illustration of directional and model-based methods.}

\begin{figure}
\centering \subfigure{
\includegraphics[width=0.4\linewidth]{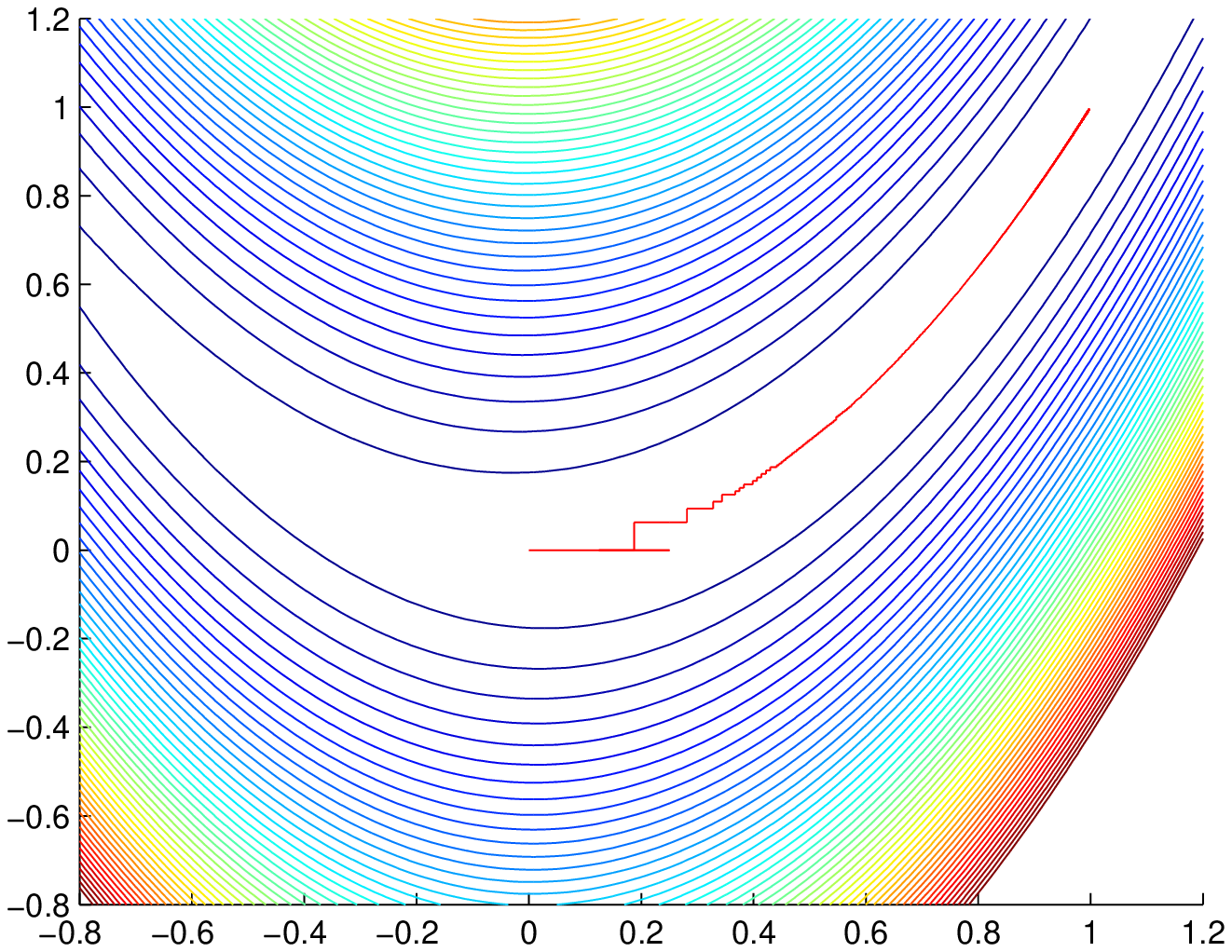}}
\centering \subfigure{
\includegraphics[width=0.4\linewidth]{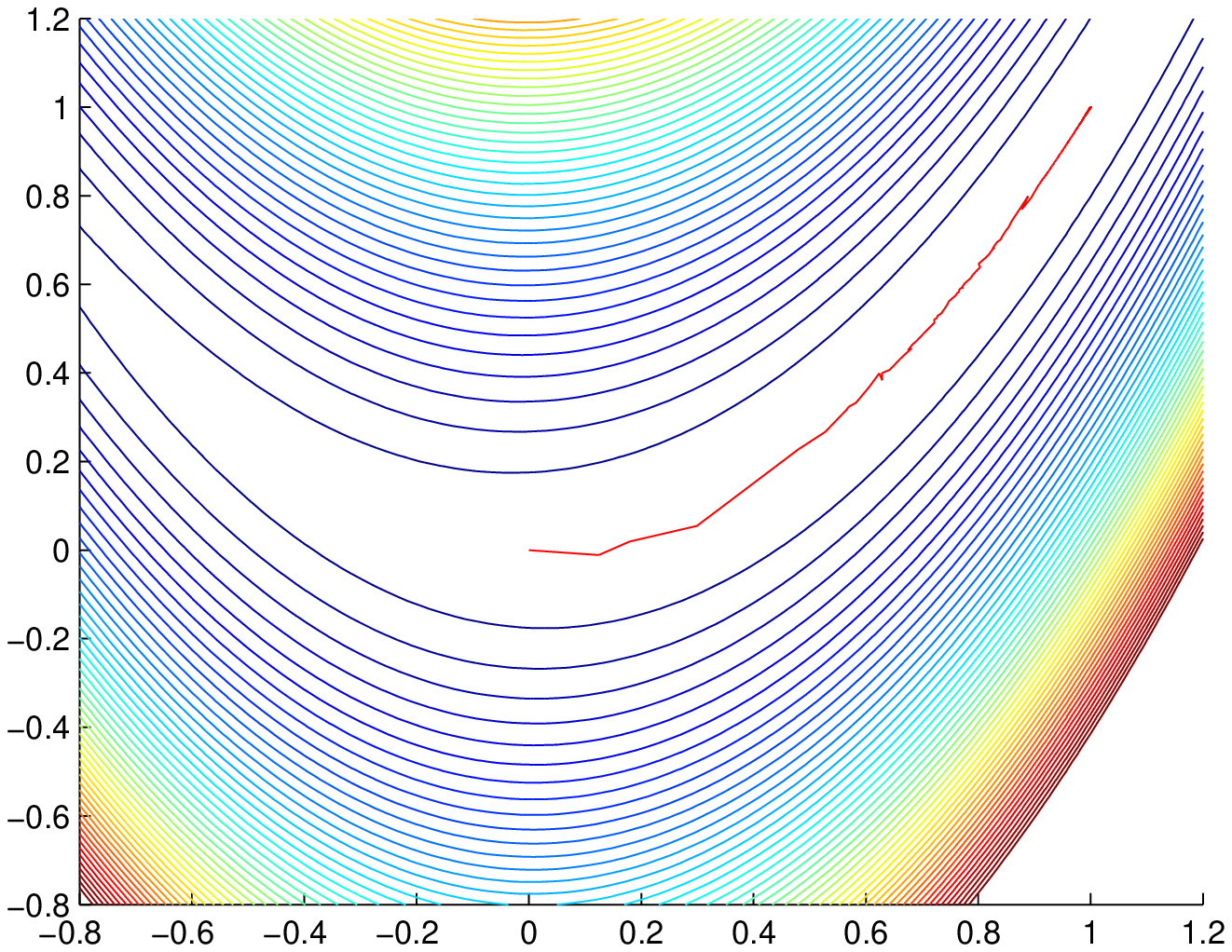}}
\caption{CS: num.eval.=$11307$, $f=10^{-6}$, DSR: num.eval.=$5756$, $f=10^{-8}$.}
\label{fig:threesearches}
\end{figure}

\begin{figure}
\centering \subfigure{
\includegraphics[width=0.4\linewidth]{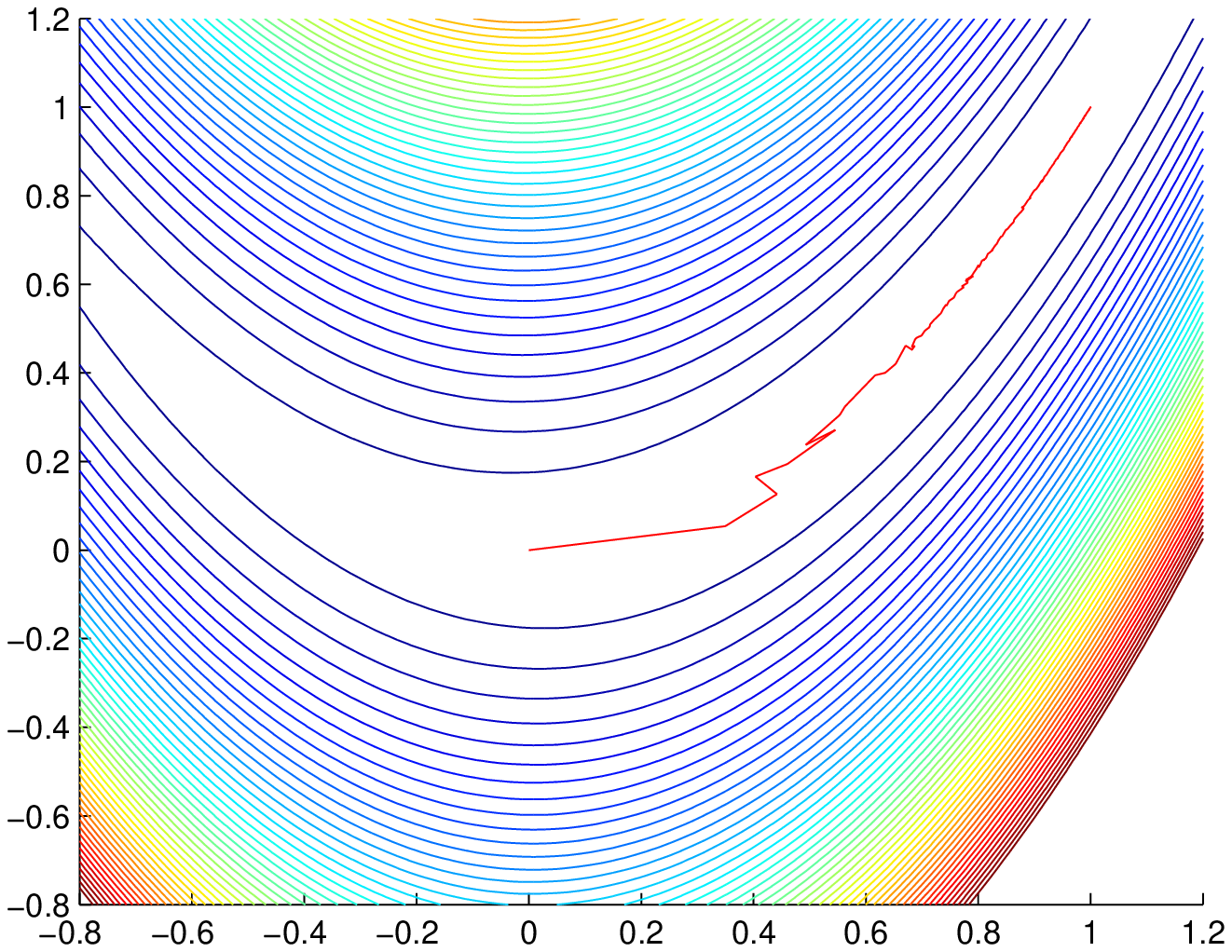}}
\centering \subfigure{
\includegraphics[width=0.4\linewidth]{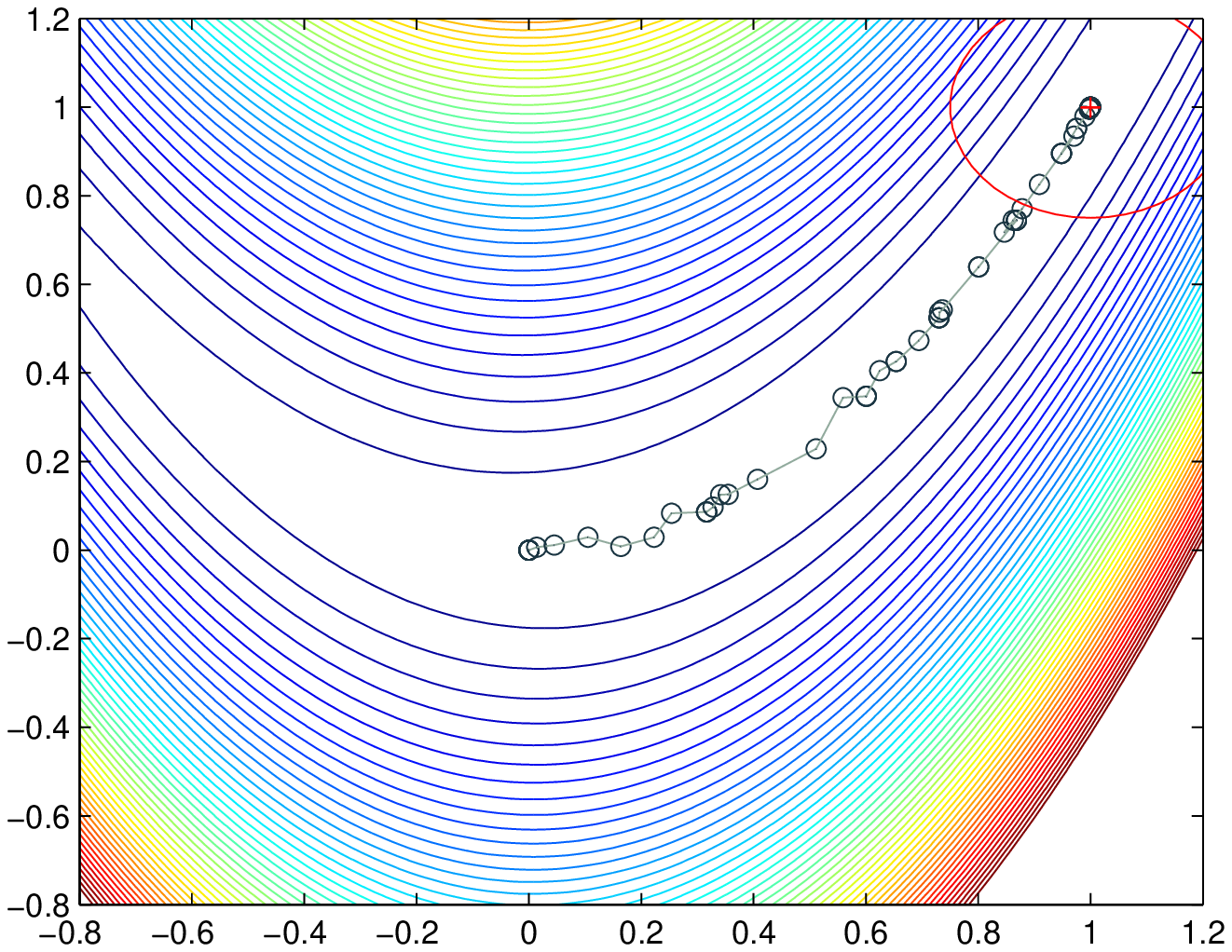}}
\caption{RS: num.eval.=$3724$, $f=10^{-8}$, TRQ: num.eval.=$62$, $f=10^{-14}$.}
\label{fig:threesearches2}
\end{figure}

Consider the well known Rosenbrock function for our computational illustration
\[
f(x) \; = \; 100(x_2 - x_1^2)^2 + (1 - x_1)^2.
\]
The function is known to be difficult for first order or zero order methods and well suited for second order methods. Nevertheless, some first/zero order methods perform reasonably, while others perform poorly.
In Figures~\ref{fig:threesearches}--\ref{fig:threesearches2} we present the contours of the function and plot the iterates produced by four methods:
1) a simple variant of direct search, the coordinate or compass search method (CS)
which uses the positive basis $[I \; -I]$,
2) a direct-search method using the positive
basis $[Q \; -Q]$ where $Q$ is an orthogonal matrix obtained by randomly generating the first column
(DSR),
3) a random search (RS) with step size inversely proportional to the iteration count,
and
4) a basic model-based trust-region method with quadratic models (TRQ). The outcome of the algorithms is summarized in the caption, which lists the number of function evaluations and the final accuracy for each method.


It is evident from these results that the random directional approaches, and in particular
random search, are more successful at finding good directions for descent, while the
coordinate search is slow due to the fixed choice of the search directions.
It is also clear, from the performance of the second order trust-region method
on this problem, that using accurate models can substantially improve efficiency.
It is natural, thus, to  consider the effects of randomization in model-based methods.
In particular we consider methods that use models built from randomly sampled
points in hopes of obtaining better models.

\section{First order trust-region method based on probabilistic models} \label{sec:first-order}

Let us consider the classical trust-region method setting and notation
(see~\cite{ARConn_KScheinberg_LNVicente_2009} for a similar description).
At iteration $k$, $f$ is   approximated by a model
$m_k$ within the ball $B(x_k,\delta_k)$ centered at $x_k$ and of radius $\delta_k$. Then
the model is minimized (or approximately
minimized)  in the ball to possibly obtain $x_{k+1}$.
In this section we will introduce and analyze a trust-region algorithm based on probabilistic models, i.e.,
 models $m_k$ which are  built in a random fashion.
First we discuss these models and state what will be assumed from them.

\subsection{The probabilistically fully linear models}

For simplicity of the presentation, we consider  only quadratic models, written in the form
\[
m_k(x_k + s ) \; = \; m_k(x_k) + s^{\top} g_k + \frac{1}{2} s^{\top} H_k s,
\]
where $g_k=\nabla m_k(x_k)$ and $H_k=\nabla^2 m_k(x_k)$.
Our analysis is not, however, dependent on the models being quadratic.

Let us start by introducing a
measure of (linear or first order)
accuracy of the model $m_k$.

\begin{definition}
We say that a function $m_k$ is a $(\kappa_{eg},\kappa_{ef})$-fully linear
model of $f$ on $B(x_k,\delta_k)$ if, for every $s\in B(0,\delta_k)$,
\begin{eqnarray*}
\|\nabla f(x_k+s) - \nabla m_k(x_k+s) \| & \leq & \kappa_{eg} \delta_k, \\
\|f(x_k+s) - m(x_k+s) \| & \leq & \kappa_{ef} \delta_k^2.
\end{eqnarray*}
\end{definition}

The concept of fully linear models is introduced in~\cite{ARConn_KScheinberg_LNVicente_2009-paper}  and~\cite{ARConn_KScheinberg_LNVicente_2009}, but here we use
the notation proposed in~\cite{SCBillups_JLarson_PGraf_2013}.
In~\cite[Chapter~6]{ARConn_KScheinberg_LNVicente_2009} there is a detailed discussion on how to construct and maintain
deterministic fully linear models.

For the case of random models, the key assumption in our convergence analysis is that these models exhibit good accuracy with sufficiently high probability. We will consider random models $M_k$, and then use the notation
$m_k= M_k(\omega_k)$ for their realizations. The randomness of the models will imply the randomness  of the current point~$x_k$ and the current trust-region radius~$\delta_k$.
Thus, in the sequel, these random quantities will be denoted by $X_k$ and $\Delta_k$, respectively,
while $x_k= X_k(\omega_k)$ and $\delta_k= \Delta_k(\omega_k)$ denote their realizations.

\begin{definition}\label{def:random_poised_mart}
We say that a sequence of random models $\{M_k\}$ is $(p)$-probabilistically
$(\kappa_{eg},\kappa_{ef})$-fully linear for a corresponding sequence $\{B(X_k,\Delta_k)\}$ if the events
\[
S_k \; = \; \{M_k \text{ is a }(\kappa_{eg},\kappa_{ef})\text{-fully linear
model of }f \text{ on }B(X_k,\Delta_k)\}
\]
satisfy the following submartingale-like condition
\[
P(S_k| F^M_{k-1}) \; \geq \; p,
\]
where $F^M_{k-1}=\sigma(M_0,\ldots,M_{k-1})$ is the $\sigma$-algebra generated by $M_0,\ldots,M_{k-1}$.
Furthermore, if $p\geq \frac12$, then we say that the random models are probabilistically $(\kappa_{eg},\kappa_{ef})$-fully linear.
\end{definition}

Note that $M_k$ is a random model that encompasses all the randomness of iteration $k$ of our algorithm. Definition~\ref{def:random_poised_mart} serves to enforce the following property:  even though the accuracy of $M_k$ may be dependent of the history of the algorithm, ($M_1,\ldots,M_{k-1}$), it is sufficiently good  with probability at least~$p$, regardless of that history. We believe this condition is more reasonable  than assuming  complete independence of $M_k$ from the past, which is difficult to ensure given that the current iterate, around which the model is built, and the trust-region radius depend on the algorithm history.

Now we discuss the corresponding assumptions on the models realizations that we use in the algorithm.
The first assumption guarantees that we are able to adequately minimize (or reduce) the model
at each iteration of our algorithm.

\begin{assumption} \label{assumptions:fcd}
For every $k$, and for all realizations $m_k$ of $M_k$ (and of $X_k$ and $\Delta_k$), we are able to compute a step $s_k$ such that
\beqn{fcd-final}
m_k(x_k) - m_k(x_k+s_k) \; \geq \; \frac{\kappa_{fcd}}{2} \| g_k \|
\min\left\{ \frac{\|g_k\|}{\| H_k \|}, \delta_k \right\}, \eeqn
for some constant $\kappa_{fcd} \in (0,1]$.
We say in this case that $s_k$ has achieved a fraction of Cauchy decrease.
\end{assumption}

The Cauchy step itself, which is the minimizer of the quadratic model within the trust region
and along the negative model gradient~$-g_k$, trivially satisfies this property
with $\kappa_{fcd}=1$.

We also assume a uniform  bound on the model Hessians:

\begin{assumption} \label{assumptions:bhm}
There exists a  positive constant $\kappa_{bhm}$, such that
 for every $k$, the Hessians $H_k$ of all realizations $m_k$ of $M_k$ satisfy
\beqn{eq:bhm}
\|H_k\| \; \leq \; \kappa_{bhm}.
\eeqn
\end{assumption}

The above assumption is introduced for convenience. While it is possible to show our results without this assumption,
it is not restrictive in the case of fully linear models. In particular, one can construct fully linear models with arbitrarily small $\|H_k\|$
using interpolation techniques.
In the case of  models that, fortuitously, have large Hessian norms, because they are not fully linear,  we can simply set the Hessian to some other matrix of a smaller norm (or zero).

\subsection{Algorithm and basic properties}

Let us consider the following simple trust-region algorithm.

\begin{algorithm}\label{alg:DFO_Random_SIMPLE1}
Fix the positive parameters $\eta_1$, $\eta_2$, $\gamma$, $\delta_{\max}$ with $\gamma >1$. At
iteration~$k$ approximate the function $f$ in $B(x_k,\delta_k)$ by
$m_k$ and then approximately minimize $m_k$ in $B(x_k,\delta_k)$, computing $s_k$ so that
it satisfies a fraction of Cauchy decrease~\eqref{fcd-final}. Let
\[
\rho_k \; = \; \frac{f(x_k) - f(x_k+s_k)}{m(x_k) - m(x_k+s_k)}.
\]
If $\rho_k \geq \eta_1$ and $\|g_k\| \geq \eta_2\delta_k$, set $x_{k+1}=x_k
+ s_k$ and $\delta_{k+1} = \min \{\gamma \delta_k, \delta_{\max}\}$. Otherwise, set $x_{k+1}=x_k$
and $\delta_{k+1} = \gamma^{-1} \delta_k$. Increase $k$ by one and repeat the iteration.
\end{algorithm}

This  is a basic trust-region algorithm, with one specific modification: the trust-region radius is always increased if sufficient function reduction is achieved, that is the step is successful, {\it and} the trust-region radius is small compared to the norm of the model gradient.
The logic behind this update follows from the fact that the step size obtained by the model minimization is typically proportional to the norm of the model gradient, hence the trust region should be of comparable size also.  Later we will show how the algorithm can be modified to allow for the trust-region radius to remain unchanged in some iterations.

Each realization of the algorithm defines a sequence of realizations for the
corresponding random variables, in particular: $m_k=M_k(\omega_k)$,
$x_k= X_k(\omega_k)$, $\delta_k=\Delta_k(\omega_k)$.

For the purpose of proving convergence of the algorithm to first order critical points, we
assume that the function $f$ and its gradient are Lipschitz
continuous in regions considered by the algorithm realizations.
To define this region we follow the process in~\cite{ARConn_KScheinberg_LNVicente_2009-paper}.
Suppose that $x_0$ (the initial iterate)
is given. Then all the subsequent  iterates
belong to the level set
\[
L(x_0) \; = \; \left\{ x \in \mathbb{R}^n: \; f(x) \leq f(x_0)
\right\}.
\]
However, the failed iterates may lie outside this set.
 In the setting considered in this paper, all potential iterates are
restricted to the region
\[
L_{enl}(x_0) \; = \; L(x_0) \cup \bigcup_{x \in L(x_0)}
B(x,\delta_{\max}) \; = \; \bigcup_{x \in L(x_0)} B(x,\delta_{\max}),
\]
where $\delta_{\max}$ is the upper bound on the size of the trust regions, as imposed by the algorithm.

\begin{assumption} \label{assumptions:gLip} Suppose $x_0$ and
$\delta_{\max}$ are given. Assume that~$f$ is continuously
differentiable in an open set containing the set $L_{enl}(x_0)$
and that $\nabla f$ is Lipschitz continuous on $L_{enl}(x_0)$ with constant $\kappa_{Lg}$.
Assume also that~$f$ is bounded from below on~$L(x_0)$.
\end{assumption}

The following lemma states that the trust-region radius converges to zero
 regardless of the realization
of the model sequence $\{ M_k \}$ made by the algorithm,
as long as the fraction of Cauchy decrease
is achieved by the step at every iteration.

\begin{lemma}\label{lemma:DFO_RandomDelta_to_0}
For every realization of Algorithm \ref{alg:DFO_Random_SIMPLE1},
$$\lim_{k\to\infty}\delta_k \; = \; 0.$$
\end{lemma}

\proof{Suppose that $\{\delta_k\}$ does not converge to zero. Then, there exists
$\epsilon>0$ such that $\#\{k:\delta_k>\epsilon\}=\infty$. Because of the
way $\delta_k$ is updated we must have
$$\#\left\{k:\delta_k>\frac{\epsilon}\gamma,\, \delta_{k+1} \geq \delta_k\right\} \; = \; \infty,$$
in other words, there must be an infinite number of iterations on which $\delta_{k+1}$ is not decreased,
 and, for these iterations we
have $\rho\geq \eta_1$ and $\|g_k\|\geq \eta_2\frac{\epsilon}\gamma $. Therefore, because \eqref{fcd-final}  and \eqref{eq:bhm} hold,
\begin{eqnarray}\label{eq:declowbnd}
f(x_k) - f(x_k+s_k) & \geq & \eta_1\left(m(x_k) - m(x_k+s_k)\right) \nonumber \\
 & \geq & \eta_1  \; \frac{\kappa_{fcd}}{2} \| g_k \|
\min \left\{ \frac{\|g_k\|}{\kappa_{bhm}}, \delta_k \right\} \\
  & \geq & \eta_1\frac{ \kappa_{fcd}}{2}\min \left\{\frac{\eta_2}{\kappa_{bhm}}, 1 \right\} \eta_2 \frac{\epsilon^2}{\gamma^2}. \nonumber
\end{eqnarray}
This means that at each iteration where $\delta_k$ is increased, $f$  is reduced by a constant. Since $f$ is bounded from below, the number of such iterations
cannot be infinite, and hence we arrived at a contradiction.}

Another result that we use in our analysis is the following fact typical of trust-region methods,
stating that, in the presence of sufficient model accuracy, a successful step will be achieved,
provided the trust-region radius is sufficiently small relatively to the size of the model
gradient.

\begin{lemma}\label{cnd-Deltaongk}
If $m_k$ is $(\kappa_{eg}, \kappa_{ef})$-fully linear on $B(x_k,\delta_k)$ and
\[
\delta_k \; \leq \; \min\left\{ \frac{\|g_k\|} { \kappa_{bhm}},
\frac{\kappa_{fcd} (1- \eta_1) \|g_k\|} { 4\kappa_{ef}} \right\},
\]
then at the $k$-th iteration $\rho_k \geq \eta_1$.
\end{lemma}

The proof can be found in~\cite[Lemma~10.6]{ARConn_KScheinberg_LNVicente_2009}.

\section{Convergence of the first order trust-region method based on probabilistic models} \label{sec:conv-1st-order}

We now assume that the  models used in the algorithm are probabilistically fully linear,
and show our first order convergence results.
First we will state an auxiliary result from the martingale literature that will be useful in our analysis.

\begin{theorem}\label{theorem:Martingales}
 Let $G_k$ be a submartingale, i.e., a sequence of random variables
 which, for every $k$, are integrable ($\EE(|G_k|) < \infty$) and
  \[
  \EE[G_k|F^G_{k-1}] \; \geq \; G_{k-1},
  \]
  where $F^G_{k-1}=\sigma(G_0,\ldots,G_{k-1})$  is the $\sigma$-algebra generated by $G_0,\ldots,G_{k-1}$
  and $\EE[G_k|F^G_{k-1}]$ denotes the conditional expectation of $G_k$ given the past history of
  events $F^G_{k-1}$.

  Assume further that $| G_k - G_{k-1} | \leq  M <\infty$, for every $k$.
   Consider the random events $C = \{\lim_{k \to \infty} G_k$ ${\rm  exists\ and\ is\ finite}\}$
and $D = \{\limsup_{k\to \infty} {G_k} = \infty\}$. Then $P (C \cup D) = 1$.
 \end{theorem}

\proof{The theorem is a simple extension of \cite[Theorem 5.3.1]{RDurret_2010}, see~\cite[Exercise~5.3.1]{RDurret_2010}.}
Roughly speaking, this results shows that a random walk with bounded increments and an upward drift either converges to a finite limit or is unbounded from above. We will apply this result to $\log \Delta_k$ which, as we show, is a random walk with an upward drift that cannot converge to a finite limit.

\subsection{The liminf-type convergence}

As it is typical in trust-region methods, we show first that a subsequence of the iterates
drive the gradient of the objective function to zero.

\begin{theorem}\label{theorem:DFO_Random_SIMPLE1_conv}
Suppose that the model sequence $\{M_k\}$ is probabilistically $(\kappa_{eg},\kappa_{ef})$-fully linear
for some positive constants $\kappa_{eg}$ and $\kappa_{ef}$.
Let $\{ X_k \}$ be a sequence of random iterates generated
by Algorithm~\ref{alg:DFO_Random_SIMPLE1}. Then, almost surely,
\[
\liminf_{k\to\infty} \|\nabla f(X_k)\| \; = \; 0.
\]
\end{theorem}

\proof{Recall the definition of events the $S_k$ in Definition \ref{def:random_poised_mart}. Let us start by constructing the following random walk
$$
W_k \; = \; \sum_{i=0}^k(2\cdot1_{S_i}-1),
$$
where $1_{S_i}$ is the indicator random variable ($1_{S_i}=1$ if $S_i$ occurs, $1_{S_i}=0$ otherwise).
From the martingale-like property enforced in Definition \ref{def:random_poised_mart}, it easily follows that $W_k$ is a submartingale.
In fact, one has
\begin{eqnarray*}
\EE[W_k|F^S_{k-1}] & = & \EE[W_{k-1}|F^S_{k-1}] + \EE[2.1_{S_k}-1|F^S_{k-1}] \\ [1ex]
& = & W_{k-1} + 2 \EE[1_{S_k}|F^S_{k-1}] - 1  \\ [1ex]
& = & W_{k-1} + 2 P ( S_k | F^S_{k-1} ) - 1  \\ [1ex]
& \geq & W_{k-1},
\end{eqnarray*}
where $F^S_{k-1}=\sigma(1_{S_0},\ldots,1_{S_{k-1}})$ is the $\sigma$-algebra generated by $1_{S_0},\ldots,1_{S_{k-1}}$,
in turn contained in
$F^M_{k-1}=\sigma(M_0,\ldots,M_{k-1})$.
Since the submartingale $W_k$ has $\pm1$ (and hence, bounded) increments it cannot have a finite limit.
 Thus, by Theorem \ref{theorem:Martingales} we have that
the event $D = \{\limsup_{k\to \infty} {W_k} = \infty\}$
holds almost surely.

Since our objective is to show that $\liminf_{k\to\infty} \|\nabla f(X_k)\| = 0$
almost surely, we can show it by conditioning on an almost sure
event. All that follows is conditioned on the event $D = \{\limsup_{k\to \infty} {W_k} = \infty\}$.

Suppose there exist $\epsilon>0$ and $k_1$ such that, with positive probability,
$$
\|\nabla f(X_k)\| \; \geq \; \epsilon,
$$
for all $k \geq k_1$.

Let  $\{ x_k \}$ and $\{ \delta_k \}$ be any realization
of $\{ X_k \}$ and $\{ \Delta_k \}$, respectively, built by Algorithm~\ref{alg:DFO_Random_SIMPLE1}.
By Lemma \ref{lemma:DFO_RandomDelta_to_0}, there exists $k_2$ such that
we have $\forall {k\geq k_2}$
\begin{equation}\label{eq:Deltamin}
\delta_k \; < \; \BBB:=\min\left\{\frac{\epsilon}{2\kappa_{eg}},\frac{\epsilon}{2\kappa_{bhm}},
\frac{\epsilon}{2\eta_2},\frac{\kappa_{fcd}(1-\eta_1)\epsilon}{8\kappa_{ef}}, \frac{\delta_{\max}}{\gamma}\right\}.
\end{equation}
Consider some iterate $k\geq k_0:=\max \{k_1,k_2 \}$ such that $1_{S_k} = 1$ (model $m_k$ is fully linear). Then, from the definition of fully linear models
\[
\|\nabla f(x_k) - g_k \| \; \leq \; \kappa_{eg} \delta_k \; < \; \frac{\epsilon}2,
\]
hence,
\[
\|g_k \| \; \geq \; \frac{\epsilon}2.
\]
Using Lemma \ref{cnd-Deltaongk}  we obtain
$\rho_k \geq \eta_1$. Also
\[
\|g_k\| \; \geq \; \frac{\epsilon}2 \; \geq \; \eta_2\delta_k.
\]
Hence, by the construction of the algorithm, and the fact that
$\delta_k\leq \frac{\delta_{\max}}{\gamma}$, we have $\delta_{k+1} = \gamma \delta_k$.

Let us consider now the random variable $R_k$ with realization $r_k = \log_{\gamma}(\BBB^{-1}\delta_k)$.
For every realization~$\{ r_k \}$ of~$\{R_k\}$ we have seen that there exists~$k_0$
such that $r_k <0$ for $k\geq k_0$. Moreover, if $1_{S_k}=1$ then
$r_{k+1} = r_k + 1$, and if $1_{S_k}=0$, $r_{k+1} \geq r_k - 1$ (implying that
$R_k$ is a submartingale). Hence,
$r_k-r_{k_0}\geq w_k-w_{k_0}$ ($w_k$ denoting a realization of $W_k$). Since we are conditioning on the event~$D$,
we have that $R_k$ has to be positive infinitely often with probability one,
contradicting the fact that for all realizations of~$\{ r_k \}$ of~$\{R_k\}$ there  exists~$k_0$
such that $r_k<0$ for $k\geq k_0$. Thus, conditioning on~$D$ we always have
that $\liminf_{k\to\infty} \|\nabla f(X_k)\| = 0$ with probability one. Therefore
\[
\liminf_{k\to\infty} \|\nabla f(X_k)\| \; = \; 0
\]
almost surely.}

\subsection{The lim-type convergence}

In this subsection we show that $\lim_{k\to\infty} \|\nabla f(X_k) \|=0$ almost surely.
Before stating and  proving the main theorem we state and prove two auxiliary lemmas.

\begin{lemma}\label{lemma:DFO_Random_Borisloses}
Let $\{Z_k\}_{k\in\NN}$ be a sequence of non-negative uniformly bounded random variables
and $\{ B_k \}$ be a sequence of
Bernoulli random variables (taking values $1$ and $-1$) such that
\[
P(B_k=1| \, \sigma(B_1,\ldots,B_{k-1}), \sigma(Z_1,\ldots,Z_k) ) \; \geq \; 1/2.
\]
Let
${\cal P}$ be the set of natural numbers $k$ such that $B_k=1$ and $\cal N=\NN \setminus {\cal P}$
(note that ${\cal P}$ and ${\cal N}$ are random sequences). Then
\[
\Prob\left(\left\{\sum_{ i\in {\cal P}} Z_i < \infty
\right\}\cap\left\{\sum_{i \in {\cal N}}Z_i=\infty\right\}\right) \; = \; 0.
\]
\end{lemma}

\proof{Let us construct the following process
$G_k= G_{k-1}+B_k Z_k$. It is easy to check that $G_k$ is a submartingale with bounded increments $\{B_k Z_k\}$. Hence we can apply Theorem~\ref{theorem:Martingales} and observe that the event $\{\limsup_{k\to \infty} {G_k} = -\infty\}$ has probability zero. Noting that $G_k= \sum_{i \in{ \cal P}, i \leq k} Z_i- \sum_{i\in{ \cal N}, i \leq k}Z_i$ and hence $\left\{\sum_{ i\in {\cal P}} Z_i < \infty \right\}\cap\left\{\sum_{i \in {\cal N}}Z_i=\infty\right\}$ implies that $\{\limsup_{k\to \infty} {G_k} = -\infty\}$, if the latter happens with probability zero, then so is the former.}

\begin{lemma}\label{auxLemma_to_theorem:DFO_Random_SIMPLE1_convLIM}
Let $\{ X_k \}$ and $\{ \Delta_k \}$ be sequences of random iterates and
random trust-region radii generated
by Algorithm~\ref{alg:DFO_Random_SIMPLE1}. Fix $\epsilon>0$ and define the sequence $\left\{K_i\right\}$ consisting of the natural numbers $k$ for which $\|\nabla f(X_k)\| >\epsilon$ (note that $K_i$ is a sequence of random variables). Then,
$$\sum_{k\in \{K_i\}}{\Delta_k} \; < \; \infty$$
almost surely.
\end{lemma}

\proof{Let $\{ m_k \}$, $\{ x_k \}$, $\{ \delta_k \}$, $\{ k_i \}$ be realizations of $\{M_k\}$, $\{X_k\}$, $\{\Delta_k\}$,
$\{ K_i \}$ respectively.
Let us separate $\{k_i\}$ in two subsequences: $\{p_i\}$ is the subsequence of $\{k_i\}$ such that $m_{p_i}$ is $(\kappa_{eg}, \kappa_{ef})$-fully linear on $B(x_{p_i}, \delta_{p_i})$, and $\{n_i\}$ is the subsequence of the remaining elements of~$\{k_i\}$.

We will now show that $\sum_{j \in \{p_i\}} {\delta_j}< \infty$ for any such realization.
If $\{p_i\}$ is finite, then this result trivially follows. Otherwise, since $\delta_k\to 0$, we have that
for sufficiently large $p_i$, $\delta_{p_i}< \BBB$, with $\BBB$ defined by \eqref{eq:Deltamin}.
Since $\|\nabla f(x_{p_i})\|>\epsilon$, and $m_{p_i}$ is fully linear on $B(x_{p_i}, \Delta_{p_i})$, then
by the derivations in Theorem~\ref{theorem:DFO_Random_SIMPLE1_conv}, we have $\|g_{p_i}\|\geq\frac {\epsilon}{2}$,
and by Lemma~\ref{cnd-Deltaongk} $\rho_{p_i} \geq \eta_1$. Hence, for all $p_i$ large enough, the decrease in the function value satisfies
$$
f(x_{p_i})-f(x_{p_i+1}) \; \geq \; \eta_1 \frac{\kappa_{fcd}}{2} \frac{\epsilon}{2} \delta_{p_i}.
$$
Thus
$$
\sum_{j \in \{p_i\}}{\delta_{j}} \; \leq \; \frac{4(f(x_0)-f_*)}{\eta_1 \kappa_{fcd}\epsilon} \; < \;\infty,
$$
where $f_*$ is a lower bound on the values of~$f$ on $L(x_0)$.

For each $k_i$, the event $S_{k_i}$ (whether the model is fully linear on iteration $k_i$) has probability at least $\frac12$ conditioned on all of the history of the algorithm.
Hence, we can apply Lemma~\ref{lemma:DFO_Random_Borisloses} (note that $\{ \Delta_k \}$ is a sequence
of non-negative uniformly bounded variables and $S_{k_i}$ are the Bernoulli random variables) and obtain
\[
\Prob\left( \left\{ \sum_{j\in \{P_i\}}{\Delta_j}< \infty \right\} \cap  \left\{
\sum_{j\in \{N_i\}}{\Delta_j}= \infty \right\} \right) \; = \;  0.
\]
This means that, almost surely,
\[
\sum_{j\in \{K_i\}}{\Delta_{j}} \; = \; \sum_{j \in \{P_i\}}{\Delta_j} +
\sum_{j \in \{N_i\}}{\Delta_j} \; < \; \infty.
\]
}

We are now ready to prove the $\lim$-type result.

\begin{theorem}\label{theorem:DFO_Random_SIMPLE1_convLIM}
Suppose that the model sequence $\{M_k\}$ is probabilistically $(\kappa_{eg},\kappa_{ef})$-fully linear
for some positive constants $\kappa_{eg}$ and $\kappa_{ef}$.
Let $\{ X_k \}$ be a sequence of random iterates generated
by Algorithm~\ref{alg:DFO_Random_SIMPLE1}. Then, almost surely,
\[
\lim_{k\to\infty} \|\nabla f(X_k)\| \; = \; 0.
\]
\end{theorem}

\proof{Suppose that $\lim_{k\to\infty} \|\nabla f(X_k)\| = 0$ does not hold almost surely. Then,
with positive probability, there exists $\epsilon>0$ such that $\|\nabla f(X_k)\|>2\epsilon$, holds for infinitely many~$k$. Without loss of generality, we assume that $\epsilon=\frac1{n_\epsilon}$, for some natural number $n_\epsilon$.

Let $\{K_i\}$ be a subsequence of the iterations for which $\|\nabla f(X_k)\|>\epsilon$.
We are going to show that, if such an $\epsilon$ exists then $\sum_{j \in \{K_i\}}{\Delta_j}$ is a divergent sum.

Let us call a pair of integers $(W^{\prime},W^{\prime\prime})$ an ``ascent'' pair if
$0<W^{\prime}<W^{\prime\prime}$, $\|\nabla f(X_{W^{\prime}})\|\leq \epsilon$, $\|\nabla
f(X_{W^{\prime}+1})\|> \epsilon$, $\|\nabla f(X_{W^{\prime\prime}})\|>2\epsilon$ and,
moreover, for any $w\in (W^{\prime},W^{\prime\prime})$, $\epsilon<\|\nabla
f(X_{w})\|\leq2\epsilon$. Each such ascent pair  forms a nonempty interval of integers $\{W^{\prime}+1, \ldots, W^{\prime\prime}\}$ which is a subset of the sequence $\{K_i\}$.  Since
$\liminf_{k\to\infty} \|\nabla f(X_k)\| = 0$ holds almost surely (by
Theorem~\ref{theorem:DFO_Random_SIMPLE1_conv}), it follows that there
are infinitely many such intervals. Let us consider the sequence of these intervals $\{(W_{\ell}^\prime, W_{\ell}^{\prime\prime})\}$. The idea is now to show (with positive probability)
that, for any ascent pair $(W_{\ell}^{\prime},W_{\ell}^{\prime\prime})$ with $\ell$ sufficiently large,
$\sum_{j=W_{\ell}^{\prime}+1}^{W_{\ell}^{\prime\prime}-1}\Delta_j$ is uniformly bounded away from $0$ (and hence $W_{\ell}^{\prime}+1<W_{\ell}^{\prime\prime}$),
which implies that   $\sum_{j \in \{K_i\}}{\Delta_j} = \infty$ since $\sum_{\ell}\sum_{j=W_{\ell}^{\prime}+1}^{W_{\ell}^{\prime\prime}-1}\Delta_j\leq
\sum_{j \in \{K_i\}}{\Delta_j}$, because the sequence $\{K_i\}$  contains all intervals $\{W_{\ell}^\prime, W_{\ell}^{\prime\prime}\}$.

Let $\{x_k\}$ and $\{ \delta_k\}$ be realizations of $\{X_k\}$ and $\{ \Delta_k\}$,
for which $\|\nabla f(x_k)\|>\epsilon$ for $k \in \{k_i\}$.
By the triangular inequality, for any $j$,
\[
\epsilon \; < \; \left| \|\nabla f(x_{w_{\ell}^{\prime}})\| - \|\nabla
f(x_{w_{\ell}^{\prime\prime}})\|\right| \; \leq \; \sum_{j=w_{\ell}^{\prime}}^{w_{\ell}^{\prime\prime}-1}\left| \|\nabla f(x_j)\|
- \|\nabla f(x_{j+1})\|\right|.
\]
Since $\nabla f$ is Lipschitz continuous (with constant $\kappa_{Lg}$),
\begin{eqnarray}
\epsilon & \leq & \sum_{j=w_{\ell}^{\prime}}^{w_{\ell}^{\prime\prime}-1}\left| \|\nabla f(x_j)\| -
\|\nabla f(x_{j+1})\|\right|\\
 &\leq& \kappa_{Lg}\sum_{j=w_{\ell}^{\prime}}^{w_{\ell}^{\prime\prime}-1}\|x_j - x_{j+1}\|\\
 &\leq& \kappa_{Lg}\left(\delta_{w_{\ell}^{\prime}}+\sum_{j=w_{\ell}^{\prime}+1}^{w_{\ell}^{\prime\prime}-1}\delta_j\right).
\end{eqnarray}
From the fact that $\delta_k$ converges to zero,
then, for any $\ell$ large enough, $\delta_{w_{\ell}^{\prime}}<\frac{\epsilon}{2\kappa_{Lg}}$,
 and hence $\sum_{j=w_{\ell}^{\prime}+1}^{w_{\ell}^{\prime\prime}-1}\delta_j >\frac{\epsilon}{2}>0$, which gives us  $\sum_{j \in \{k_i\}}{\delta_j} = \infty$.

We have thus proved that if, $\lim_{k\to\infty} \|\nabla f(X_k)\| = 0$ does not hold almost surely, then, with positive probability, there exists $n_\epsilon$ such that $\{K_i\}$ defined as above based on $n_\epsilon$,  satisfies $\sum_{j \in \{K_i\}}{\Delta_j} = \infty$.

On the other hand, Lemma~\ref{auxLemma_to_theorem:DFO_Random_SIMPLE1_convLIM} guarantees  that, for every $n_\epsilon$, the probability of $\sum_{j \in \{K_i\}}{\Delta_j} = \infty$ is zero. Since there are countable $n_\epsilon\in \mathbb{N}$ and since the union of a countable number of rare events is still rare we have that the probability of existence of a value $n_\epsilon$ for which $\sum_{j \in \{K_i\}}{\Delta_j} = \infty$ is zero, which contradicts the initial assumption that $\lim_{k\to\infty} \|\nabla f(X_k)\| = 0$ does not hold almost surely.}

\subsection{Modified trust-region schemes}

The trust-region radius update of Algorithm~\ref{alg:DFO_Random_SIMPLE1} may be too restrictive as it only allows for this radius to be increased or decreased. In practice typically two separate thresholds are used, one for the increase of the trust-region radius and another for its decrease. In the remaining cases the trust-region radius remains unchanged.  Hence, here we propose an algorithm similar to Algorithm~\ref{alg:DFO_Random_SIMPLE1} but slightly more appealing in practice.

\begin{algorithm}
Fix the positive parameters $\eta_1$, $\eta_2$, $\eta_3$, $\gamma$, $\delta_{\max}$, with $\gamma >1$ and $\eta_3\leq \eta_2$. At
iteration~$k$ approximate the function $f$ in $B(x_k,\delta_k)$ by
$m_k$ and then approximately minimize $m_k$ in $B(x_k,\delta_k)$, computing $s_k$ so that
it satisfies a fraction of Cauchy decrease~\eqref{fcd-final}. Let
\[
\rho_k \; = \; \frac{f(x_k) - f(x_k+s_k)}{m(x_k) - m(x_k+s_k)}.
\]
If $\rho_k \geq \eta_1$, then set $x_{k+1}=x_k + s_k$ and
\[
\delta_{k+1} \; = \;
\left
\{\begin{array}{ll}
          \gamma^{-1} \delta_k & \text{if} \;\; \|g_k\| < \eta_3 \delta_k,   \\
\\
          \delta_k & \text{if} \;\; \eta_3 \delta_k \leq \|g_k\| < \eta_2 \delta_k,   \\
\\
          \min\{\gamma \delta_k, \delta_{\max}\} & \text{if} \;\; \eta_2 \delta_k \leq \|g_k\|.
        \end{array}
\right.
\]
Otherwise, if $\rho_k \leq \eta_1$, set $x_{k+1}=x_k$ and $\delta_{k+1} = \gamma^{-1} \delta_k$.
\end{algorithm}

It is straightforward to adapt the proofs of Lemma \ref{lemma:DFO_RandomDelta_to_0} and Theorems \ref{theorem:DFO_Random_SIMPLE1_conv} and \ref{theorem:DFO_Random_SIMPLE1_convLIM} to show
the  convergence for this new algorithm. Additionally, one can consider two different thresholds, $0<\eta_0<1$ for decrease of the trust region radius, and $\eta_1>\eta_0$ for the increase of the trust region radius.

\section{Second order trust-region method based on probabilistic models} \label{sec:2nd-order}

In this section we present the analysis of the convergence of a trust-region algorithm to second order stationary points  under the assumption that the random models are likely to provide second order accuracy.

\subsection{The probabilistically fully quadratic  models}

Let us now introduce a measure of second order quality or accuracy of the models $m_k$
(see~~\cite{ARConn_KScheinberg_LNVicente_2009-paper, ARConn_KScheinberg_LNVicente_2009,SCBillups_JLarson_PGraf_2013} for more details).

\begin{definition} \label{def:FQ}
We say that a function $m_k$ is a $(\kappa_{eh},\kappa_{eg},\kappa_{ef})$-fully quadratic
model of $f$ on $B(x_k,\delta_k)$ if, for every $s\in B(0,\delta_k)$,
\begin{eqnarray*}
\|\nabla^2 f(x_k+s) - H_k \| & \leq & \kappa_{eh} \delta_k, \\
\|\nabla f(x_k+s) - \nabla m_k(x_k+s) \| & \leq & \kappa_{eg} \delta_k^2, \\
\|f(x_k+s) - m(x_k+s) \| & \leq & \kappa_{ef} \delta_k^3.
\end{eqnarray*}
\end{definition}

As in the fully linear case, we assume that the models used in the algorithms
are fully quadratic with a certain probability.

\begin{definition}\label{def:random_poised_mart_quad}
We say that a sequence of random models $\{M_k\}$ is $(p)$-probabilistically
$(\kappa_{eh},\kappa_{eg},\kappa_{ef})$-fully quadratic for a corresponding sequence $\{B(X_k,\Delta_k)\}$ if the events
\[
S_k \; = \; \{M_k \text{ is a }(\kappa_{eh},\kappa_{eg},\kappa_{ef})\text{-fully quadratic
model of }f \text{ on }B(X_k,\Delta_k)\}
\]
satisfy the following submartingale-like condition
\[
P(S_k| F^M_{k-1}) \; \geq \; p,
\]
where $F^M_{k-1}=\sigma(M_0,\ldots,M_{k-1})$ is the $\sigma$-algebra generated by $M_0,\ldots,M_{k-1}$.
Furthermore, if $p\geq \frac12$, then we say that the random models are probabilistically $(\kappa_{eh},\kappa_{eg},\kappa_{ef})$-fully
quadratic.
\end{definition}

We now need to discuss the algorithmic requirements and problem assumptions
which will be needed for global convergence to second order critical points.
In terms of problems assumptions we will need one more order of smoothness.

\begin{assumption} \label{assumptions:Lip2_cont}
Suppose $x_0$ and $\delta_{\max}$ are given.
Assume that $f$ is
twice continuously differentiable in an open set containing the
set $L_{enl}(x_0)$ and that $\nabla^2 f$ is Lipschitz continuous with constant $\kappa_{Lh}$ and that
$\|\nabla^2 f\|$ is bounded by a constant $\kappa_{bhf}$ on
$L_{enl}(x_0)$.
Assume also that~$f$ is bounded from below on~$L(x_0)$.
\end{assumption}

We will no longer assume that the Hessian $H_k$ of the models is bounded in norm, since we cannot simply disregard
large Hessian model values without possibly affecting the chances of the model being fully quadratic. However, a simple analysis can show that $\|H_k\|$ is uniformly bounded from above for any
fully quadratic model $m_k$ (although we may not know what this bound is and hence may not be able to use it in an algorithm).

\begin{lemma}\label{lemma:bmh}
Given constants $\kappa_{eg}$, $\kappa_{ef}$, $\kappa_{eh}$, and $\delta_{\max}$,
there exists a constant $\kappa_{bmh}$ such that
for every~$k$ and every realization~$m_k$ of $M_k$
which is a $(\kappa_{eg},\kappa_{ef}, \kappa_{eh})$-fully quadratic
model of $f$ on $B(x_k,\delta_k)$ with $x_k\in L(x_0)$ and $\delta_k\leq \delta_{\max}$ we have
\[
\|H_k\| \; \leq \; \kappa_{bmh}.
\]
\end{lemma}

\proof{The proof follows trivially from the definition of fully quadratic models and the assumption that
$\|\nabla^2 f\|$ is bounded by a constant $\kappa_{bhf} $ on
$L_{enl}(x_0)$.}

It will also be necessary to assume that the minimization of the model
achieves a certain level of second order improvement (an extension of the Cauchy decrease).

\begin{assumption} \label{assumptions:fod}
For every $k$, and for all realizations $m_k$ of $M_k$ (and of $X_k$ and $\Delta_k$), we are
able to compute a step $s_k$ so that
\beqn{fod-final} m_k(x_k) - m_k(x_k+s_k) \; \geq \;
\frac{\kappa_{fod}}{2} \max \left\{
            \|g_k\|
            \min\left[ \frac{\|g_k\|}{\| H_k \|}, \delta_k \right],
           \max \{-\lambda_{\min}(H_k), 0\} \delta_k^2
           \right\}.
 \eeqn for some constant $\kappa_{fod}
\in (0,1]$.
We say in this case that $s_k$ has achieved a fraction of optimal decrease.
\end{assumption}

A step satisfying this assumption is given, for instance, by
computing both the Cauchy step and,
in the presence of negative curvature in the model,
the eigenstep, and by choosing the
one that provides the largest reduction in the model.
The eigenstep is the minimizer of the quadratic model in the trust region along an eigenvector
corresponding to the smallest (negative) eigenvalue of~$H_k$.

The measure of proximity to a second order stationary point
for the function~$f$ is slightly different from the traditional,
and is given by
\[
\tau(x) \;=\;
\max \left\{ \min \left[ \| \nabla f(x) \|, \frac{\| \nabla f(x) \|} {\| \nabla^2 f(x) \|} \right]
, -\lambda_{\min}(\nabla^2f(x)) \right\}.
\]
The model approximation of this measure is defined similarly:
\[
\tau^m(x) \;=\;
\max \left\{ \min \left[ \| \nabla m(x) \|, \frac{\| \nabla m(x) \|} {\| \nabla^2 m(x) \|}\right],
-\lambda_{\min}(\nabla^2m(x)) \right\}.
\]
We consider the additional terms $\| \nabla f(x) \| / \| \nabla^2 f(x) \|$ and
$\| \nabla m(x) \| / \| \nabla^2 m(x) \|$ given that we no longer assume a bound
in the model Hessians as we did in the first order case.
We show now that $\tau(x)$ is Lipschitz continuous under Assumption~\ref{assumptions:Lip2_cont}.

\begin{lemma}\label{lem:tau_Lip}
Suppose that
Assumption~\ref{assumptions:Lip2_cont} holds.
Then there exists a constant $\kappa_{L\tau}$ such that for all $x_1, x_2\in L_{enl}(x_0)$
\beqn{tau_Lip}
|\tau(x_1) - \tau(x_2)| \; \leq \; \kappa_{L\tau} \|x_1-x_2\|.
\eeqn
\end{lemma}

\proof{First we note that under Assumption~\ref{assumptions:Lip2_cont} there must  exist
an upper bound $\kappa_{bfg} > 0$ on the norm of the gradient of~$f$,
$\| \nabla f(x) \| \leq \kappa_{bfg}$ for all $x \in L_{enl}(x_0)$.

Then let us see that
$h(x) = \min \{ \| \nabla f(x) \|, \| \nabla f(x) \| / \| \nabla^2 f(x) \| \}$
is Lipschitz continuous. Given $x,y \in L_{enl}(x_0)$, one consider four cases:
(i) The case $\| \nabla^2 f(x) \| \geq 1$ and $\| \nabla^2 f(y) \| \geq 1$
results from the Lipschitz continuity and boundedness above of the gradient and the Hessian.
(ii) The case $\| \nabla^2 f(x) \| < 1$ and $\| \nabla^2 f(y) \| < 1$
results from the Lipschitz continuity of the gradient.
(iii) The argument is the same for the other two cases, so let us choose one of
them, say $\| \nabla^2 f(x) \| < 1$ and $\| \nabla^2 f(y) \| \geq 1$. In this case,
using these inequalities, one has
\begin{eqnarray*}
|h(x)-h(y)| \; \leq \; \left|
\|\nabla f(x)\| - \frac{\| \nabla f(y) \|} {\| \nabla^2 f(y) \|} \right|
\; \leq \;
\left|
\|\nabla f(x)\| - \frac{\| \nabla f(x) \|} {\| \nabla^2 f(y) \|} \right| + \frac{\kappa_{Lg} \|x-y\|}{\| \nabla^2 f(y) \|} \\
 \leq \;
\| \nabla f(x) \| (\| \nabla^2 f(y) \| - \| \nabla^2 f(x) \|)+ \kappa_{Lg} \|x-y\|.
\end{eqnarray*}
Thus,
$|h(x)-h(y)| \leq (\kappa_{bfg} \kappa_{Lh} + \kappa_{Lg} ) \|x-y\|$.

The proof then results  from the fact the maximum of two Lipschitz
continuous functions is Lipschitz
continuous and the fact that eigenvalues are Lipschitz continuous
functions of the entries of a matrix.}

The following lemma shows that the
difference between the problem measure $\tau(x)$ and the model
measure $\tau^m(x)$ is of the order of $\delta$ if $m(x)$ is a fully quadratic model on $B(x, \delta)$
(thus extending the error bound on the Hessians given in Definition~\ref{def:FQ}).

\begin{lemma}\label{tau_diff}
Suppose that
Assumption~\ref{assumptions:Lip2_cont} holds.
Given constants  $\kappa_{eg}$, $\kappa_{ef}$, $\kappa_{eh}$, and $\delta_{\max}$ there exists a constant $\kappa_{e\tau}$  such that for any $m_k$ which is $(\kappa_{eg},\kappa_{ef}, \kappa_{eh})$-fully quadratic
model of $f$ on $B(x_k,\delta_k)$ with $x_k\in L(x_0)$ and $\delta_k\leq \delta_{\max}$ we have
\beqn{tau_error} |\tau(x_k) - \tau^m(x_k)| \; \leq \;
\kappa_{e\tau} \delta_k. \eeqn
\end{lemma}

\proof{From the definition of fully quadratic models and the upper bounds
on $\|\nabla f\|$ and $\|\nabla^2 f\|$ on $L_{enl}(x_0)$, we conclude that both $\|\nabla m (x_k)\|$
and $\|\nabla^2 m (x_k)\|$ are also bounded from above with constants independent of $x_k$ and $\delta_k$.

For a given $x_k$ several situation may occur depending on which terms dominate in the expressions for
$\tau(x_k)$ and $\tau^m(x_k)$.  In particular, if  $\| \nabla^2 f(x_k) \|\leq 1$ and $\| \nabla^2 m(x_k) \|\leq 1$,
then
\[
\tau(x_k) \;=\;
\max \left \{\| \nabla f(x_k) \|,
 -\lambda_{\min}(\nabla^2f(x_k)) \right\}
\]
and
\[
\tau^m(x_k) \;=\;
\max \left \{ \| \nabla m(x_k) \|,
-\lambda_{\min}(\nabla^2m(x_k)) \right\}
\]
and   the proof of the lemma is the same as in the case of the usual criticality measure, analyzed in~\cite{ARConn_KScheinberg_LNVicente_2009}.
Let us consider the case, when $\| \nabla^2 f(x_k) \|\geq 1$ and $\| \nabla^2 m(x_k) \|\geq 1$.
From the fact that $m_k$ is $(\kappa_{eg},\kappa_{ef}, \kappa_{eh})$-fully quadratic we have that
\begin{eqnarray*}
\left| \frac{\| \nabla m(x_k) \|} {\| \nabla^2 m(x_k) \|}-  \frac{\| \nabla f(x_k) \|} {\| \nabla^2 f(x_k) \|}\right|
 & \leq & \left|\| \nabla m(x_k) \| \| \nabla^2 f(x_k) \|-  \| \nabla f(x_k) \| \| \nabla^2 m(x_k) \| \right|
 \\ & \leq& \| \nabla^2 f(x_k) \| \kappa_{ef}  \delta^2_k
+  \| \nabla f(x_k) \| \kappa_{eh}\delta_k \;\; \leq \;\;  \kappa_{e\tau} \delta_k,
 \end{eqnarray*}
for some large enough $\kappa_{e\tau}$, independent of $x_k$ and $\delta_k$.

The other two cases that need consideration are
\begin{itemize}
\item $\tau^m(x_k)=  \frac{\| \nabla m(x) \|} {\| \nabla^2 m(x_k) \|}$,  $\tau(x_k)=  \| \nabla f(x_k) \|$, and
\item  $\tau(x_k)=  \frac{\| \nabla f(x) \|} {\| \nabla^2 f(x_k) \|}$,  $\tau^m(x_k)=  \| \nabla m(x_k) \|$.
\end{itemize}
Let us consider the first case
\begin{eqnarray*}
|\tau(x_k)-\tau^m(x_k)| \;=\; \left| \|\nabla f(x_k)\|- \frac{\|\nabla m(x_k)\|}{\| \nabla^2 m(x_k)\|}
\right| \;\leq\; \left|   \|\nabla f(x_k)\|- \frac{\|\nabla f(x_k)\|}{\| \nabla^2 m(x_k)\|}
\right|  +\frac{\kappa_{eg}\delta^2_k}{\| \nabla^2 m(x_k)\|} \;\leq \\
\|\nabla f(x_k)\| (\| \nabla^2 m(x_k)\|-1+\kappa_{eg}\delta^2_k ) \;\leq\;  \|\nabla f(x_k)\| (\kappa_{eh}\delta_k+\kappa_{eg}\delta^2_k) \;\leq\; \kappa_{e\tau} \delta_k,
 \end{eqnarray*}
for some large enough $\kappa_{e\tau}$, independent of $x_k$ and $\delta_k$.  The proof of the second case is derived in a similar manner. Combining these results with standard steps of analysis, such as
the one in  in~\cite{ARConn_KScheinberg_LNVicente_2009} we conclude the proof of this lemma.}

Let us now define $\tau_k=\tau(x_k)$ and $\tau^m_k=\tau^m(x_k)$.
From the assumption that $\nabla^2 f(x)$ is bounded on $L_{enl}(x_0)$, it is clear
that if  $\tau_k \to 0$ (when $k \to \infty$), then $\nabla f(x_k)\to  0$ and
$\max\{-\lambda_{\min}(\nabla^2 f(x_k)), 0\}\to 0$.
We next present an algorithm for which we will then analyze the convergence of $\tau_k$.

\subsection{Algorithm and liminf-type convergence}

Consider the following modification of Algorithm~\ref{alg:DFO_Random_SIMPLE1}.

\begin{algorithm}\label{alg:DFO_Random_SIMPLE2}
Fix the positive parameters $\eta_1$, $\eta_2$, $\gamma$, with $\gamma >1$. At
iteration~$k$ approximate the function $f$ in $B(x_k,\delta_k)$ with
$m_k$ and then approximately minimize $m_k$ in $B(x_k,\delta_k)$, computing $s_k$ so that
it satisfies a fraction of optimal decrease~\eqref{fod-final}. Let
\[
\rho_k \; = \; \frac{f(x_k) - f(x_k+s_k)}{m(x_k) - m(x_k+s_k)}.
\]
If $\rho_k \geq \eta_1$ and $\tau_k \geq \eta_2\delta_k$, set $x_{k+1}=x_k
+ s_k$ and $\delta_{k+1} = \min \{\gamma \delta_k, \delta_{\max}\}$. Otherwise, set $x_{k+1}=x_k$
and $\delta_{k+1} = \gamma^{-1} \delta_k$. Increase $k$ by one and repeat the iteration.
\end{algorithm}

The analysis of this method is similar to that of the first order method described in Section~\ref{sec:first-order}.
The main difference lies in a replacement of the use of assumptions and in the lack of proof of the $\lim$-type result.
First, we will follow the steps of Section~\ref{sec:first-order} to analyze the behavior
of the trust-region radius.

\begin{lemma}\label{lemma:DFO_RandomDelta_to_0_2}
 For every realization of Algorithm \ref{alg:DFO_Random_SIMPLE2},
$$\lim_{k\to\infty}\delta_k \;=\; 0.$$
\end{lemma}

\proof{Suppose that $\{\delta_k\}$ does not converge to zero. Then, there exists
$\epsilon>0$ such that $\#\{k:\delta_k>\epsilon\}=\infty$. We are going to consider the following subsequence
$\{k:\delta_{k}>\frac{\epsilon}\gamma,\,
\delta_{k+1} \geq \delta_k\}$. By assumption this subsequence is infinite and due to the
way $\delta_k$ is updated we have  $\tau^m_k \geq \eta_2\frac{\epsilon}\gamma$  for each $k$ in this subsequence.

First assume that $\min \{\| g_k \| ,\frac{ \| g_k\|}{\|H_k\|}\} \geq \eta_2\frac{\epsilon}\gamma$.
Therefore, from \eqref{fod-final} we have
\begin{eqnarray*}
f(x_k) - f(x_k+s_k) & \geq & \eta_1\left(m(x_k) - m(x_k+s_k)\right) \\
 & \geq & \eta_1 \frac{\kappa_{fod}}{2} \|g_k\| \min \left\{ \frac{\|g_k\|}{\| H_k \|},\delta_k \right\} \\
 & \geq & \eta_1 \frac{\kappa_{fod}}{2} \eta_2\frac{\epsilon^2}{\gamma^2} \min\{\eta_2,1\} .
\end{eqnarray*}

Now assume that $ -\lambda_{\min}(H_k) \geq \eta_2\frac{\epsilon}\gamma$.
Therefore, from \eqref{fod-final} we have
\begin{eqnarray*}
f(x_k) - f(x_k+s_k) & \geq & \eta_1\left(m(x_k) - m(x_k+s_k)\right) \\
 & \geq & - \eta_1 \frac{\kappa_{fod}}{2}   \lambda_{\min}(H_k) \delta_k^2 \\
 & \geq & \eta_1 \frac{\kappa_{fod}}{2}  \eta_2 \frac{\epsilon^3}{\gamma^3}.
\end{eqnarray*}

This means that at iteration $k$ the function $f$ decreases by an amount bounded
away from zero. Since we have assumed that there is an infinite number of
such iterations, we obtain a contradiction.}

The next step is to extend Lemma~\ref{cnd-Deltaongk} to the second order context.

\begin{lemma}\label{2cnd-Deltaontauk}
If $m_k$ is $(\kappa_{eh}, \kappa_{eg}, \kappa_{ef})$-fully quadratic on $B(x_k,\delta_k)$ and
\[ \delta_{k}
\; \leq \; \min\left \{ \tau_{k}^m,
\sqrt{\frac{\kappa_{fod} ( 1 - \eta_1 ) \tau_k^m} {
4\kappa_{ef}}}, \frac{ \kappa_{fod} ( 1 -
\eta_1 ) \tau_k^m} { 4\kappa_{ef}}\right \}, \]
then at the $k$-th iteration $\rho_k\geq \eta_1$.
\end{lemma}

The proof is a trivial extension of the proof of~\cite[Lemma 10.17]{ARConn_KScheinberg_LNVicente_2009} taking into account our modified definition
of $\tau_{k}^m$ .

We can now prove the following convergence result which states that a subsequence of iterates
approaches second order stationarity almost surely.

\begin{theorem}\label{theorem:DFO_Random_SIMPLE2_conv}
Suppose that the model sequence $\{M_k\}$ is probabilistically $(\kappa_{eh},\kappa_{eg},\kappa_{ef})$-fully quadratic
for some positive constants  $\kappa_{eh}$, $\kappa_{eg}$, and $\kappa_{ef}$.
Let $\{ X_k \}$ be a sequence of random iterates generated
by Algorithm~\ref{alg:DFO_Random_SIMPLE2}. Then, almost surely,
\[
\liminf_{k\to\infty} \tau_k \; = \; 0.
\]
\end{theorem}

\proof{As in Theorem~\ref{theorem:DFO_Random_SIMPLE1_conv},
let us consider the random walk $W_k = \sum_{i=0}^k(2\cdot1_{S_i}-1)$
(where $1_{S_i}$ is the indicator random variable, now based on the event $S_i$ of Definition~\ref{def:random_poised_mart_quad}).
All that follows is also conditioned on the almost sure event $D = \{\limsup_{k\to \infty} {W_k} = \infty\}$.

Suppose there exist $\epsilon>0$ and $k_1$ such that, with positive probability,
$\tau_k \geq \epsilon$, for all $k \geq k_1$.
Let $\{ x_k \}$ and $\{ \delta_k \}$ be any realization
of $\{ X_k \}$ and $\{ \Delta_k \}$, respectively, built by Algorithm~\ref{alg:DFO_Random_SIMPLE2}.
From Lemma \ref{theorem:DFO_Random_SIMPLE2_conv}, there exists $k_2$ such that
we have $\forall_{k\geq k_2}$
\begin{equation} \label{b-2nd}
\Delta_k \; < \; \BBB :=\min\left\{\frac{\epsilon}{2\kappa_{e \tau}},\frac{\epsilon}{2},\frac{\epsilon}{2\eta_2},
\sqrt{\frac{\kappa_{fod}  ( 1 - \eta_1 ) \epsilon} {
8\kappa_{ef}}}, \frac{ \kappa_{fod} ( 1 -
\eta_1 ) \epsilon} { 8\kappa_{ef}},  \frac{\delta_{\max}}{\gamma} \right\} \;> \; 0.
\end{equation}
Let $k\geq k_0:=\max \{k_1,k_2 \}$ such that $1_{S_k} = 1$. Then,
$|\tau_k - \tau^m_k | \leq \kappa_{e \tau} \delta_k < \frac{\epsilon}2$, and thus
$\tau^m_k \geq \frac{\epsilon}2$.
Now, using Lemma~\ref{2cnd-Deltaontauk}, we obtain
$\rho_k \geq \eta_1$. We also have $\tau^m_k \geq \frac{\epsilon}2 \geq \eta_2\delta_k$.
Hence, by the construction of the algorithm, and the fact that
$\delta_k\leq \frac{\delta_{\max}}{\gamma}$, we have $\delta_{k+1} = \gamma \delta_k$.

The rest of the proof is derived exactly as the proof of Theorem~\ref{theorem:DFO_Random_SIMPLE1_conv}
(defining the random variable~$R_k$ with realization $r_k = \log_{\gamma}(\BBB^{-1}\delta_k)$,
but with  $\BBB$ now given by~(\ref{b-2nd})).
Conditioning on~$D$ we obtain $\liminf_{k\to\infty} \tau_k = 0$, and thus
$\liminf_{k\to\infty} \tau_k = 0$ almost surely.}

\subsection{The $\lim$-type convergence}

Let us summarize what we know about the convergence of Algorithm
\ref{alg:DFO_Random_SIMPLE2}. Clearly all results that hold for Algorithm
\ref{alg:DFO_Random_SIMPLE1} also hold for Algorithm
\ref{alg:DFO_Random_SIMPLE2}, hence as long as the probabilistically fully linear (or fully quadratic) models are used, almost surely,
 the iterates  of Algorithm \ref{alg:DFO_Random_SIMPLE2} form a sequence $\{x_k\}$, such that $\nabla f(x_k)\to 0$ as $k\to \infty$, in other words,
the sequence $\{x_k\}$ converges to a set of first order stationary points. Moreover, as we just showed in the previous section,
as long as the probabilistically  fully quadratic models are used, there exists a subsequence of iterates  $\{x_k\}$ which converges to a second order stationary point with probability one.
Note that under certain assumptions, for instance, assuming that the Hessian of $f(x)$
is strictly positive definite at every second order stationary point, we can conclude
from the results shown so far
(and similarly to \cite[Theorem 6.6.7]{ARConn_NIMGould_PhLToint_2000})
that, almost surely, all limit points of the sequence of iterates of
Algorithm~\ref{alg:DFO_Random_SIMPLE2} are second order stationary points.

There are however cases, when the set of first order stationary points is connected, and contains both second order stationary points and points with negative curvature of the Hessian.
An example of such a function is
\[
f(x)=xy^2.
\]
All points such that $y=0$ for a set of first order stationary points, while any $x\geq 0$ gives us second order stationary points, while  $x< 0$ does not.
In theory our algorithm may produce two subsequences of iterates, one converging to a point with $y=0$ and $x>0$ (a second order stationary point), and another converging to a point for which $y=0$ and  $x<0$ (a first order stationary point with negative curvature of the Hessian).

Theorem 6.6.8 in~\cite{ARConn_NIMGould_PhLToint_2000} shows that all
limit points of a trust-region algorithm are second order stationary
without the assumption on these limit points being isolated, but under the
condition that the trust-region radius is increased at successful
iterations. The results in~\cite{ARConn_KScheinberg_LNVicente_2009} show
that all limit points of a
 trust-region framework based on deterministic fully quadratic models are
second order stationary under a slightly modified trust-region maintenance
conditions.  While the same result may be true for
Algorithm~\ref{alg:DFO_Random_SIMPLE2} using probabilistically fully quadratic models, we were unable to
extend the results in~\cite{ARConn_KScheinberg_LNVicente_2009}  to this case.
Below we present explanations where such extension fails, but the key lies
in the fact that successful iterations and hence increase in the trust
region are no longer guaranteed.

\begin{conjecture}\label{theorem:DFO_Random_SIMPLE2_convLIM}
Suppose that the model sequence $\{M_k\}$ is probabilistically $(\kappa_{eh}, \kappa_{eg},\kappa_{ef})$-fully quadratic
for some positive constants $\kappa_{eh}$, $\kappa_{eg}$ and $\kappa_{ef}$.
Let $\{ X_k \}$ be a sequence of random iterates generated
by Algorithm~\ref{alg:DFO_Random_SIMPLE2}. Then, almost surely,
\[
\lim_{k\to\infty} \tau_k \; = \; 0.
\]
\end{conjecture}

Let us attempt to follow the same logic as in the proof of Theorem \ref{theorem:DFO_Random_SIMPLE1_convLIM}.  The first part of the proof applies immediately after substituting
$\|\nabla f(x)\|$ by $\tau(x)$ wherever is appropriate.

Indeed, suppose that $\lim_{k\to\infty} \tau(X_k) = 0$ does not hold almost surely. Then,
with positive probability, there exists $\epsilon>0$ such that $\tau(X_k)>2\epsilon$, holds for infinitely many~$k$'s.
 Without loss of generality, we assume that $\epsilon=\frac1{n_\epsilon}$, for some natural number $n_\epsilon$.

Let $\{K_i\}$ be a subsequence of the iterations for which $\tau(X_k)>\epsilon$.
We are going to show that, if such an $\epsilon$ exists then $\sum_{j \in \{K_i\}}{\Delta_j}$ is a divergent sum.

Let us call a pair of integers $(W^{\prime},W^{\prime\prime})$ an ``ascent'' pair if
$0<W^{\prime}<W^{\prime\prime}$, $\tau(X_{W^{\prime}})\leq \epsilon$, $\tau(X_{W^{\prime}+1})> \epsilon$,
$\tau (X_{W^{\prime\prime}})>2\epsilon$ and,
moreover, for any $w\in (W^{\prime},W^{\prime\prime})$, $\epsilon<\tau(X_{w})\leq2\epsilon$.
Each such ascent pair  forms a nonempty interval of integers $\{W^{\prime}+1, \ldots, W^{\prime\prime}\}$ which is a subset of the sequence $\{K_i\}$.  Since
$\liminf_{k\to\infty} \tau(X_k) = 0$ holds almost surely (by
Theorem~\ref{theorem:DFO_Random_SIMPLE2_conv}), it follows that there
are infinitely many such intervals. Let us consider the sequence of these intervals $\{(W_{\ell}^\prime, W_{\ell}^{\prime\prime})\}$. The idea is now to show (with positive probability)
that, for any ascent pair $(W_{\ell}^{\prime},W_{\ell}^{\prime\prime})$ with $\ell$ sufficiently large,
$\sum_{j=W_{\ell}^{\prime}+1}^{W_{\ell}^{\prime\prime}-1}\Delta_j$ is uniformly bounded away from $0$ (and hence $W_{\ell}^{\prime}+1<W_{\ell}^{\prime\prime}$),
which implies that   $\sum_{j \in \{K_i\}}{\Delta_j} = \infty$ since $\sum_{\ell}\sum_{j=W_{\ell}^{\prime}+1}^{W_{\ell}^{\prime\prime}-1}\Delta_j\leq
\sum_{j \in \{K_i\}}{\Delta_j}$, because the sequence $\{K_i\}$  contains all intervals $\{W_{\ell}^\prime, W_{\ell}^{\prime\prime}\}$.

Let $\{x_k\}$ and $\{ \delta_k\}$ be realizations of $\{X_k\}$ and $\{ \Delta_k\}$,
for which $\tau_k > \epsilon$ for $k \in \{k_i\}$.
By the triangular inequality, for any $j$,
\[
\epsilon \; < \; \left| \tau_{w_{\ell}^{\prime}} - \tau_{w_{\ell}^{\prime\prime}}\right| \; \leq \; \sum_{j=w_{\ell}^{\prime}}^{w_{\ell}^{\prime\prime}-1}\left| \tau_j
- \tau_{j+1} \right|.
\]
Since $\tau(x)$ is Lipschitz continuous (with constant $\kappa_{L\tau}$),
\begin{eqnarray}
\epsilon & \leq & \sum_{j=w_{\ell}^{\prime}}^{w_{\ell}^{\prime\prime}-1}\left| \tau_j-
\tau_{j+1} \right|\\
 &\leq& \kappa_{L\tau}\sum_{j=w_{\ell}^{\prime}}^{w_{\ell}^{\prime\prime}-1}\|x_j - x_{j+1}\|\\
 &\leq& \kappa_{L\tau}\left(\delta_{w_{\ell}^{\prime}}+\sum_{j=w_{\ell}^{\prime}+1}^{w_{\ell}^{\prime\prime}-1}\delta_j\right).
\end{eqnarray}
From the fact that $\delta_k$ converges to zero,
then, for any $\ell$ large enough, $\delta_{w_{\ell}^{\prime}}<\frac{\epsilon}{2\kappa_{L\tau}}$,
 and hence $\sum_{j=w_{\ell}^{\prime}+1}^{w_{\ell}^{\prime\prime}-1}\delta_j >\frac{\epsilon}{2}>0$, which gives us  $\sum_{j \in \{k_i\}}{\delta_j} = \infty$.

We have thus proved that if, $\lim_{k\to\infty} \tau(X_k) = 0$ does not hold almost surely, then, with positive probability, there exists $n_\epsilon$ such that $\{K_i\}$ defined as above based on $n_\epsilon$,  satisfies $\sum_{j \in \{K_i\}}{\Delta_j} = \infty$.

The second part of the proof should rely on providing the contradiction to the statement that $\sum_{j \in \{K_i\}}{\Delta_j} = \infty$ can happen with positive probability. In the case of Theorem~\ref{theorem:DFO_Random_SIMPLE1_convLIM} the proof utilized  the fact that the sum of all $\delta_j$, for all $j \in  \{p_i\}$ such that $m_j$ is probabilistically fully linear, has to be finite, because
it appears as a term in the lower bound on the total decrease of the objective function (see Lemma \ref{auxLemma_to_theorem:DFO_Random_SIMPLE1_convLIM}). However, in the second order
case the total decrease of the objective function is bounded from below by a factor of $\sum_{j \in \{p_i\}} \delta_j^2$. Hence it is possible that  $\sum_{j \in \{k_i\}} \delta_j^2<\infty$, while
 $\sum_{j \in \{k_i\}}{\delta_j}=\infty$. In the deterministic case the proof of the fact that $\sum_{j \in \{k_i\}}{\delta_j}<\infty$ relies on the trust-region maintenance strategy. In particular if the model $m_j$ is
 fully quadratic for every $j\in \{k_i\}$ then the trust-region radius is increased at each iteration  $j\in \{k_i\}$. In other words, for large enough $\ell$,  $\delta_{w_{\ell}^{\prime}+2}=\gamma \Delta_{w_{\ell}^{\prime}+1}$
 and so on until $\delta_{w_{\ell}^{\prime\prime}-1}=\gamma \delta_{w_{\ell}^{\prime\prime}-2}$. Hence
   \begin{equation}\label{eq:sequencebound}
  \sum_{j=w_{\ell}^{\prime}+1}^{w_{\ell}^{\prime\prime}-1} {\delta_j} \; < \;  M \delta_{w_{\ell}^{\prime\prime}-1}.
  \end{equation}
This then implies that  $ \sum_{j=w_{\ell}^{\prime}+1}^{w_{\ell}^{\prime\prime}-1} {\delta_j}\to 0$ as $\ell\to 0$ because $ \delta_{w_{\ell}^{\prime\prime}-1}\to 0$.

The difficulty in the probabilistic case comes from the fact that the trust region only increases with some high probability, but not necessarily at each iteration. In general it is possible to construct examples
of  random walks, which satisfy the conditions on the maintenance of
$\delta_j$, but for which with positive probability there exist arbitrarily large indices $w_{\ell}^{\prime}$ and $w_{\ell}^{\prime\prime}$, such that  $w_{\ell}^{\prime\prime}-w_{\ell}^{\prime}$  is arbitrarily large, and such that (\ref{eq:sequencebound}) does not hold. In other words, under the current conditions we are not able to show that (\ref{eq:sequencebound})  holds with probability one for all $\ell$, hence we cannot prove the conjecture.  The proof can be established under the additional assumption that the probability of occurrence of fully quadratic models increases in  some cases as the algorithm progresses. However, this scenario essentially leads to deterministic schemes and hence is of no interest for this paper.

\section{Examples of probabilistically fully linear and fully quadratic models in  DFO}\label{sec:rand_models}

In the previous section we described an algorithmic framework that is based on models whose approximation quality is random and
is sufficiently good with probability more than $1/2$ (conditioned on the past). We called such models probabilistically fully linear or fully quadratic, depending on
the quality of approximation that they provide. In this section, we discuss how such models can be generated (for some large enough values of the $\kappa$ constants) and
outline future  research in this direction.

\subsection{Fully linear  and fully quadratic polynomial interpolation models and $\Lambda$-poised sample sets}

Let  $\calP_n^d$ denote the set of polynomials of degree $\leq d$ in
$\R^n$ and let $q_1=q+1$ denote the dimension of this space. It clear
that the dimension of $\calP_n^1$ is $q_1=n+1$ and the dimension of
$\calP_n^2$ is $q_1=\half (n+1)(n+2)$.  A basis $\Phi=
\{ \phi_0(x), \phi_1(x), \ldots, \phi_{q}(x) \}$ for $\calP_n^d$ is a set of
$q_1$ polynomials of degree $\leq d$ that span $\calP_n^d$. For any such basis
$\Phi$, any polynomial $m(x)\in \calP_n^d$ can be written as
\begin{equation}
\label{eq:modelviaalpha}
m(x) \; = \; \sum_{j=0}^{q}\alpha_j\phi_j(x),
\end{equation}
where the $\alpha_j$'s are real coefficients.
Given a set of $p_1=p+1$ points $Y=\{y_0, y_1, \ldots, y_{p}\}\subset
\Re^n$, $m(x)$ is said to be the interpolation polynomial of $f(x)$ on $Y$ if it satisfies
\begin{equation}\label{eq:mbphiY}
M(\Phi, Y) \alpha \; = \; f(Y),
\end{equation}
where $M(\Phi, Y)$ is defined as follows
\begin{equation}\label{MphiY}
M(\Phi, Y) \; = \; \left [ \begin{array}{cccc}\phi_0(y^0) & \phi_1(y^0) & \cdots & \phi_{q}(y^0) \\
\phi_0(y^1) & \phi_1(y^1) & \cdots & \phi_{q}(y^1) \\
\vdots & \vdots & \vdots & \vdots\\
\phi_0(y^{p}) & \phi_1(y^{p}) & \cdots & \phi_{q}(y^{p}) \end{array}\right ]
\end{equation}
and $f(Y)$ is the $p_1$ dimensional vector whose entries are $f(y_i)$ for $i=0, \ldots, p$.
The interpolation polynomial $m(x)$ exists and
is unique if and only if $p=q$ and the set $Y$ is
\emph{poised}~\cite{MJDPowell_2001}, which essentially means that $M(\Phi, Y)$ is nonsingular.
When the number of points $p_1$ is smaller than the number of elements in $\Phi$ the matrix
$M(\Phi, Y)$  has more columns than rows and the system \eqref{eq:mbphiY} is underdetermined.
In this case there are several choices of interpolating polynomials, which we will discuss later.
If, on the other hand, $p>q$, then the system  \eqref{eq:mbphiY} is overdetermined and one can apply least squares
  regression instead of interpolation. Other polynomial approximations are also possible.
If  $Y$ is  such  that the condition number of  $M(\Phi, \hat Y)$ is  bounded by $\Lambda$, where  $\hat Y=\{(y_0-x_k)/\Delta, \ldots, (y_p-x_k)/\Delta)$ is a
 scaled version of $Y$, then we say that $Y$  is $\Lambda$-poised (see~\cite{ARConn_KScheinberg_LNVicente_2009}).
It is shown in \cite{ARConn_KScheinberg_LNVicente_2008a,ARConn_KScheinberg_LNVicente_2008b,ARConn_KScheinberg_LNVicente_2009} that if $Y$ is $\Lambda$-poised and $p_1\geq n+1$, then
one can build a   model which is $(\kappa_{ef}$, $\kappa_{eg})$-fully linear, with $\kappa_{ef}$ and $\kappa_{eg}$ both equal to $\mathcal{O}(p\Lambda)$. Analogously, it is shown that
if  $p_1\geq (n+1)(n+2)/2$, then
one can also build a $(\kappa_{ef}$, $\kappa_{eg}$, $\kappa_{eh})$-fully quadratic polynomial model  with $\kappa_{ef}$, $\kappa_{eg}$, and $\kappa_{eh}$  equal to $\mathcal{O}(p\Lambda)$.
Hence, to build
a fully quadratic model in $n$ dimension one may require
$(n+1)(n+2)/2$ sample points
(within reasonable proximity of the current iterate $x_k$).
If such a sample set is already available,
estimating the condition number
of the matrix $M(\Phi, Y)$ may require up to $\mathcal{O}(n^6)$ arithmetic operations.
This dependency on the dimension limits the use of  fully quadratic models
to  small dimensional problems.

There are two main ways to improve the per-iteration complexity of DFO algorithms. One approach, to only change the sample set by one point at a time,
has been very successful in practice, as it not only reduces the number of function evaluations, but also the linear algebra involved  \cite{ARConn_KScheinberg_PhLToint_1998,GFasano_JLMorales_JNocedal_2009,MJDPowell_2004,SMWild_2008}.
 However, in \cite{KScheinberg_PhLToint_2010} is was shown that such algorithms still require computing a $\Lambda$-poised set in the criticality step of the trust-region framework.
 Hence  computation of~$n$ new sample points is required if fully linear  models are used, while for fully quadratic models   $(n+1)(n+1)/2-1$ new sample points have to be evaluated.

The other, complementary, approach is to use quadratic models based on fewer than $(n+1)(n+2)/2$ sample points, which also reduces both  the cost of the linear algebra and the number of function evaluations.  In practical DFO applications, incomplete quadratic models have been
used very successfully.

Let  $Y=\{y_0, y_1, \ldots,y_{p}\}$
be a set  of $p+1$ sample points  with $p<q$ and let
$\Phi = \{ 1, x_1, x_2, \ldots, x_n,$
$x_1^2/2, x_1x_2, \ldots, x_{n-1}x_n  x_n^2/2 \}$,
The interpolating polynomial for $f$ on the set $Y$ is given by \eqref{eq:modelviaalpha},
where  $\alpha$ satisfies the undetermined interpolation system \eqref{eq:mbphiY}.
Since this system admits multiple solutions we have some freedom in
selecting $\alpha$.
In \cite{ARConn_KScheinberg_PhLToint_1998,ARConn_KScheinberg_LNVicente_2009,SMWild_2008} the minimum Frobenius norm (MFN) models are considered, i.e., the models
for which the Frobenius norm of the Hessian, or
$\|\alpha_Q\|_2=\|(\alpha_{n+1}, \alpha_{n+2}, \ldots, \alpha_{q})\|_2$, is minimized subject to \eqref{eq:mbphiY}.
In \cite{MJDPowell_2004} Powell selects the model based on minimizing the
Frobenius norm of the update of the Hessian.
Both these methods are successful in practice, and
provide useful second order information.
However, so far
 theoretically they are not shown to be superior to
simple linear models. Indeed as we will show in the example below the MFN models may
be nearly as bad as simple linear models, but the use of random sample sets
can provide a significant practical improvement in this case.

\subsection{Random sample sets}
In the cases when function evaluations are not very expensive or  can be obtained in parallel, there is less incentive to reuse old sample points for model building, because ensuring the  model quality can become the bottleneck of the computations.  Instead, one can simply use well-poised deterministic sample sets, chosen in advance. However this is  not always the best approach, because the pattern is chosen without any consideration for the shape of the function and may be a very poor fit. In Section \ref{sec:dfo} we have seen examples where random directions provided better decrease on average than those from a fixed pattern.   Similarly, random sample sets can automatically provide good quality models with high enough probability, yet they do not suffer from the worst case behavior of the deterministic sample sets. We consider another example.

\paragraph*{Example of comparison of using underdetermined quadratic models based on random and deterministic sample sets.}

Let us consider the  function
$f(x)=10(x_2 - x_1^2)^2 + (1 - x_1)^2$, which is a version of the Rosenbrock function used in Section \ref{sec:dfo}, but with smaller curvature and Hessian condition number.  We apply a trust-region method~\cite{ASBandeira_KScheinberg_LNVicente_2012} to this function with models based on $5$ points at each iteration and construct MFN models based on the sample sets.  Note that a fully quadratic model requires $6$ points. In one case we choose deterministic models with the sample set selected as the current iterate plus the coordinate steps of length $\delta$, i.e.,
$Y   = \{y^0, y^0\pm \delta e_1, y^0\pm \delta e_2\}$ --- a very well poised set.
 In other words, the set $Y$ is generated around the current iterate $y^0$ by adding coordinate steps of size $\delta$.
 For the second method we generate the set $Y$ by picking $4$ random points in a ball of radius $\delta$ around the current iterate. The results are as follows:  the method based on  deterministic sample sets achieved the final function value of $10^{-4}$ in $8500$ function evaluations, while the method based on random sample sets achieved the function value of $10^{-6}$ in
  $2700$   function evaluations. Clearly, using random sample sets enhances the performance of the MFN models here.
In particular, one can observe the slow progress of the deterministic method in Table \ref{tab:dfo_deter} which represents iteration output. It is clear from the table that iterations follow a pattern
(which starts at around iteration $1000$) where $\delta_k$ increased and decreased according to alternating successful and unsuccessful steps, while the progress is slow overall.

\begin{table}\label{tab:dfo_deter}
\centerline{\begin{tabular}{|c|c|c|c|c|}
\hline
Iter. \# & success & $f$ value&  $\Delta$ & $\rho$ \\
\hline
1687 &     0   & +3.67420711e-04 & +3.12e-02 & -1.66e+00 \\
1688 &     1   & +3.67418778e-04 & +6.25e-02 & +8.13e+02 \\
1689 &     0   & +3.67418778e-04 & +3.12e-02 & -1.66e+00 \\
1690 &     1   & +3.67409693e-04 & +6.25e-02 & +3.92e+03 \\
1691 &     0   & +3.67409693e-04 & +3.12e-02 & -1.66e+00 \\
1692 &     1   & +3.67407812e-04 & +6.25e-02 & +8.34e+02 \\
1693 &     0   & +3.67407812e-04 & +3.12e-02 & -1.66e+00 \\
1694 &     1   & +3.67398959e-04 & +6.25e-02 & +4.03e+03 \\
1695 &     0   & +3.67398959e-04 & +3.12e-02 & -1.66e+00 \\ \hline
\end{tabular}
}
 \end{table}

\paragraph*{Analysis of poisedness of random sample sets.}

Let us consider a sample set $Y=\{0, y^1, \ldots, y^{p}\}\subset \R^n$ with a fixed point at the origin and the remaining~$n$ points being generated randomly from a standard Gaussian distribution centered at the origin. Let us consider
$\Phi(x)=\{1, x_1, x_2, \ldots, x_n\}$.
Hence $M(\Phi, Y)$ is simply a matrix whose first column is all $1$'s, the first row is zero except the first element and the remaining $p\times n$ matrix is a Gaussian random matrix.
Under a simple transformation, the condition number of  $M(\Phi, Y)$ is equal to the condition number of a random Gaussian $p\times n$ matrix. From results in random matrix theory \cite{ZChen_JJDongarra_2005,AEdelman_1988} we have the following bound
\[
P( \cond(M(\Phi, Y)) \; > \; \Lambda) \; \leq \; C(n, p) \frac{1}{\Lambda^{|n-p|+1}},
\]
where $C(n,p)$ is a constant dependent on $p$ and $n$.  In particular, for $p=n$ the result in \cite{ZChen_JJDongarra_2005} implies
 \[
 P( \cond(M(\Phi, Y))>\Lambda) \; \leq \; \frac{1}{\sqrt{2\pi}} \frac{Cn}{\Lambda}.
 \]
where $C$ is a universal constant smaller than $6.5$.
This result implies that given $n$ and $p$, there exists~$\Lambda$ large enough such that $P( \cond(M(\Phi, Y)<\Lambda) >\frac{1}{2}$.
Hence there exist constants $\kappa_{ef}$ and $\kappa_{eg}$ such that the linear interpolation (or regression) polynomials based on Gaussian sample sets   are probabilistically $(\kappa_{ef}$, $\kappa_{eg})$-fully linear.


A more complicated, but important class of models are the quadratic models based on sample sets of $p+1>n+1$ sample points. In this case the
basis $\Phi$ is constructed from first and second order polynomials and  ${M(\Phi, Y)}$  no longer has the simple structure   of a Gaussian matrix. Matrices of this form have been studied in \cite{HRauhut_2010} and are referred to as structured random matrices. The bounds derived in  \cite{HRauhut_2010} show that the condition number of  ${M(\Phi, Y)}$ is small with sufficiently high probability if  $p$ is large enough (but still scales pseudo-linearly with the number of columns in ${M(\Phi, Y)}$).  We  believe that these results can be used to show that for a fixed~$p$, the condition number of     ${M(\Phi, Y)}$ is bounded with sufficiently high probability. Explicit derivations of such bounds and conditions is subject for future research.

\paragraph*{Sparse models based on random sample sets.}

A natural question  that arises in our context is
 whether we can build accurate, i.e., fully quadratic models,  without
requiring $(n+1)(n+2)/2$ sample points.
For instance, in larger dimensional cases it often happens that the
Hessian of the  objective function  is sparse.
Clearly, if we know in advance that some elements of the Hessian
(coefficients of $\alpha_Q$) are zero, then
we can reduce the number of variables in system~(\ref{eq:mbphiY}).
However, in a typical situation of a black-box optimization, the information about the sparsity of the Hessian is not
available.
It has been shown recently in \cite{ASBandeira_KScheinberg_LNVicente_2012} that  by minimizing $\|\alpha_Q\|_1$ instead of $\|\alpha_Q\|_2$ it is possible to recover fully quadratic  interpolation  models of a function with sparse Hessian
 by using fewer sample points than $(n+1)(n+2)/2$. This is the first result that shows that fully quadratic model recovery with incomplete sample sets is possible.
This result relies on the theory of sparse recovery in compressed sensing \cite{ECandes_2006} and on results in random matrix theory~\cite{HRauhut_2010}.
 In particular these models are shown to be fully quadratic with probability larger than $1-n^{-\gamma\log p}$, for some universal constant $\gamma>0$, as long as the number of sample points  satisfies $p \geq \mathcal{O}\left(n(\log n)^4\right)$ and a sparse fully quadratic model exists.

 Similar, but much simpler results can be obtained for  recovery of a sparse fully linear model, if such a model exists. In this case, the sample set $Y$ can be generated by a Gaussian distribution around the current iterate and the random matrix ${M(\Phi, Y)}$ can be viewed as a Gaussian matrix, just as it described above. Sparse signal recovery can be applied
 in this well-known case to show that if the number of nonzeros in the gradient is $s$ and the number of sample points is
 \[
 p \; \geq \; Cs \log(n/s),
\]
then a sparse fully linear model can be recovered with probability greater than $1-c_1e^{-c_2 p}$, for some universal constants $c_1$, $c_2$, and $C$. In fact the constants also depend on
the error between the function values $f(y_i)$ and the sparse model values $m(y_i)$, but we omit these details here for simplicity.

\paragraph{Nonuniform recovery and martingale property.}

In the examples we considered so far the sample sets are generated to provide high quality of the models independently of the past history of the algorithm. However, our theory allows the probability of a good model to be dependent on the past. In some cases taking this into account may provide a more efficient approach to building models. Here we discuss one possible example.

The results of recovery of sparse models which we considered so far from compressed sensing imply, the so-called, uniform recovery, where the matrix  ${M(\Phi, Y)}$ is designed in such a way that {\sl any} sparse model can be recovered. However, in our case, it is sufficient to recover the specific model that happens to approximate the objective function $f$ sufficiently
well in a trust region. Thus, the {\em nonuniform} recovery results can apply. Some of these results, including the ones for the Gaussian matrices, can be found in~\cite{UAyaz_HRauhut_2013,HRauhut_2010}.
The key is that if only one fixed signal needs to be recovered with high probability, then it is sufficient to generate the random matrix  ${M(\Phi, Y)}$ using fewer samples than what is necessary for the uniform recovery. The probability of generating a fully linear or fully quadratic model can be made sufficiently high, conditioned on
the model itself. This fact, in our setting, means that the probability of a  ``good'' model is high conditioned on the current iterate and trust-region radius, in other words, on the past behavior of the algorithm. In short, we observe that such a setting will satisfy the submartingale property, but not complete independence on the past.

\paragraph*{Example of comparing performance of sparse model recovery vs. other underdetermined second order models.}

\begin{figure}
\centering \subfigure{
\includegraphics[width=0.3\linewidth]{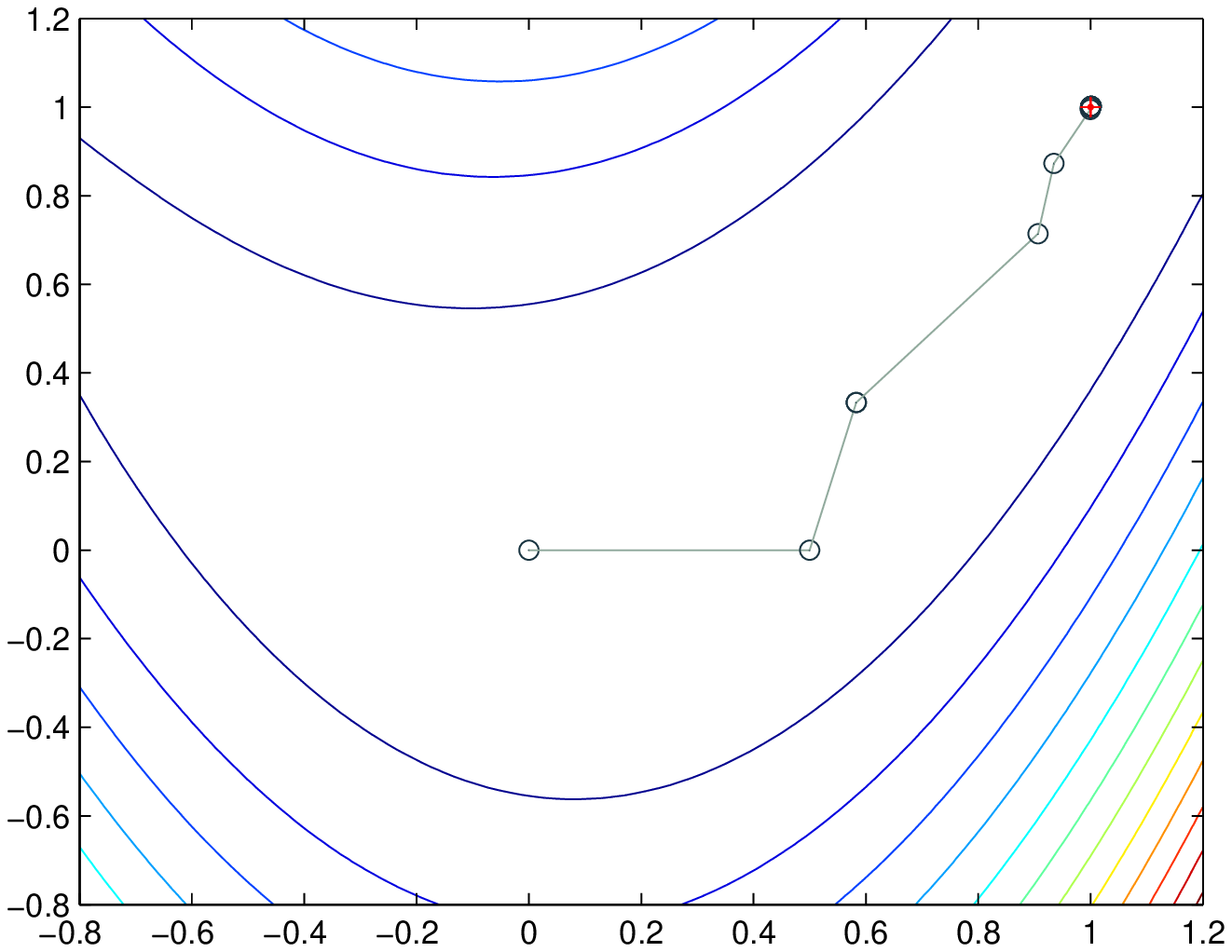}}
\centering \subfigure{
\includegraphics[width=0.3\linewidth]{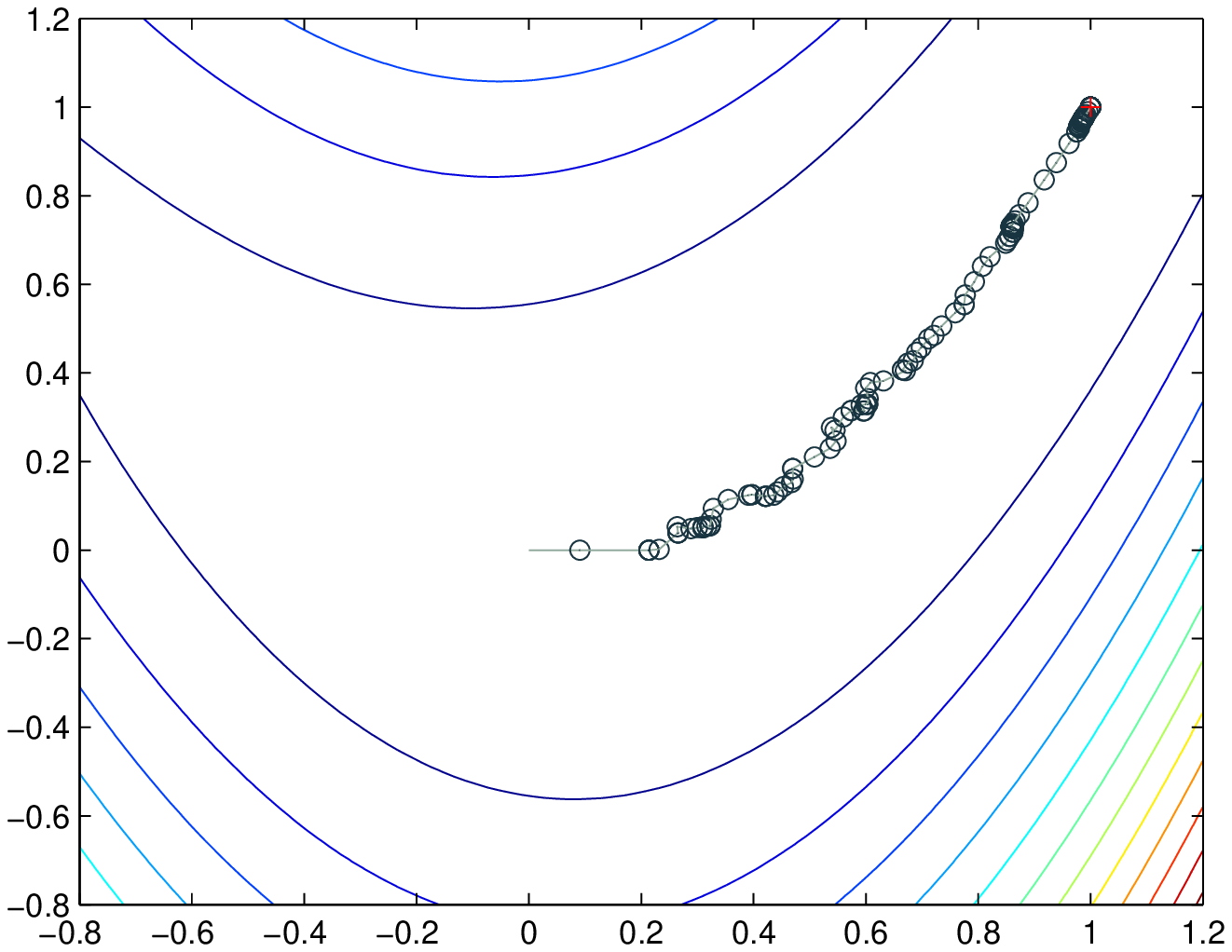}}
\centering \subfigure{
\includegraphics[width=0.3\linewidth]{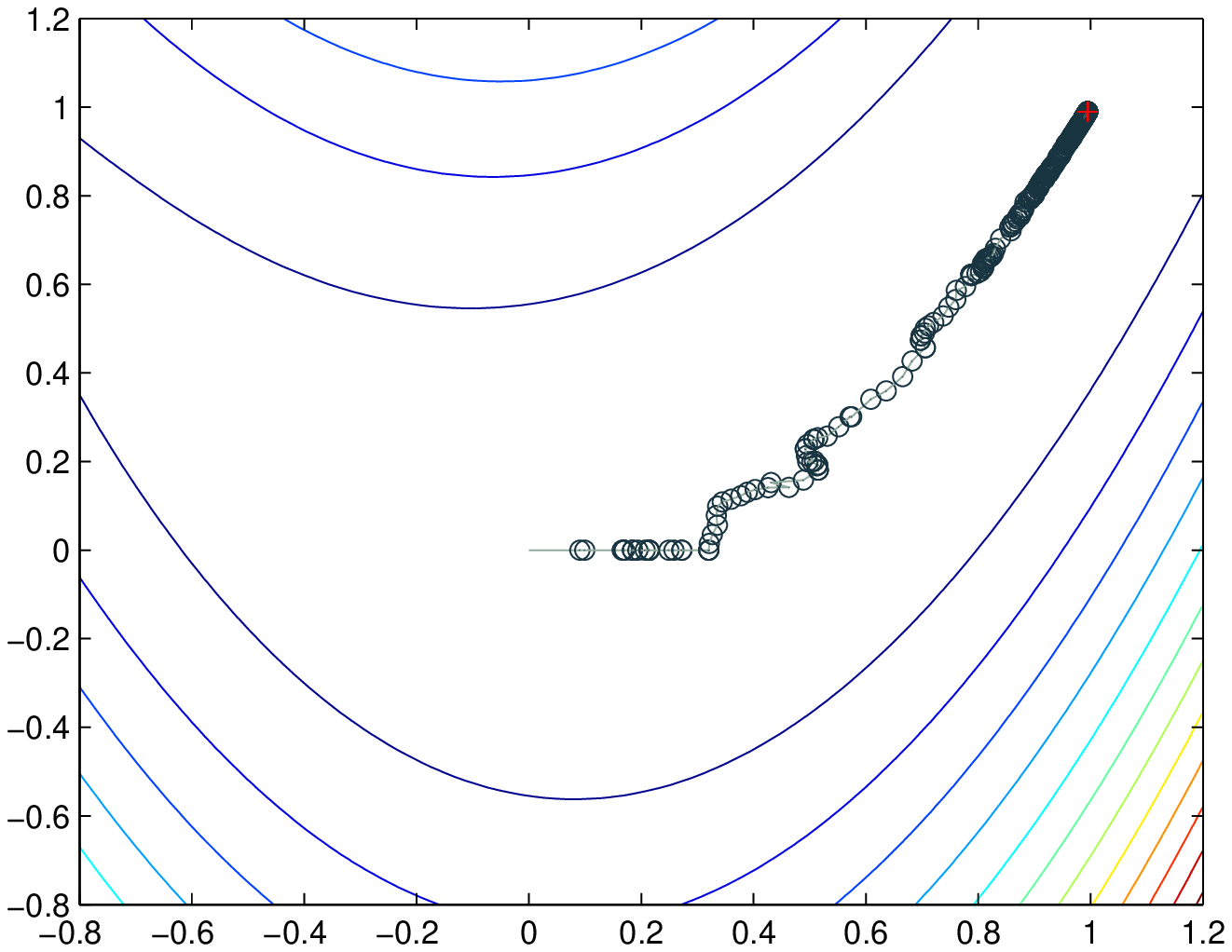}}
\caption{RSTR: {\em it}=18, {\em nf}=494, {\em f}=$10^{-11}$,  GSTR: {\em it}=164, {\em nf}=185, {\em f=}$10^{-8}$, MFN: {\em it}=325, {\em nf}=346, {\em f=}$10^{-5}$,  }
\label{fig:threedfos}
\end{figure}

Consider the following function again, $f(x)=10(x_2 - x_1^2)^2 + (1 - x_1)^2$,
but this time  $x\in \R^{10}$, which means that we have a $10$-dimensional problem, but only the first  two dimensions are important.
 Note that to build a fully linear  model
without applying sparse recovery we need to sample $11$ points and to build a fully quadratic model we need $66$ sample points.
We apply three variants of the trust-region algorithm~\cite{ASBandeira_KScheinberg_LNVicente_2012} to this problem which only differ by the choice of the models. In the first case the models are built based on $26$ random points that are
distributed in a small hypercube around the current iterate (a range of points from 20 to 30 was tried with similar results), we call this method RSTR. In the second case
we build sparse models based on ``greedy'' sample sets of up to 31 points, which only use points generated in the course of the trust-region steps, in other words, reusing  old points,
we call it GSTR. The third algorithm uses the same greedy sample sets, but constructs MFN models, rather than sparse models, we call this method MFN. The resulting optimization paths are illustrated in Figure \ref{fig:threedfos} and the final outcome is  as follows:
\begin{enumerate}
\item RSTR: Number of iterations: 18, number of function evaluations: 494, final function value:  4.0e-11.
\item GSTR: Number of iterations: 164,  number of function evaluations: 185, final function value:  5.0e-8.
\item  MFN: Number of iterations: 325,  number of function evaluations:  346, final function value:  2.5e-5.
\end{enumerate}

This example illustrates that the RSTR clearly recovers the fully quadratic models of $f$, while the other two methods do not. This is evident from the number of iterations required by each algorithm. While the first algorithm performs more function evaluations, they can be obtained in parallel, and the achieved accuracy is by far better than that of the other methods.

\paragraph*{Other random models.}

Additional settings where relying on random models may give an advantage for an optimization scheme occur in a parallel environment when full synchronization is not needed. In other words, if
function evaluations are obtained in parallel for a collection of sample points, some of the function evaluations may take much longer than others. In that case  it is possible to compute a model based on a sufficiently large subset of sampled values and ignore the points whose function values are not returned on time. Under the assumption that the function computation failures  occur randomly, the remaining subset  is still a random sample set.

Alternatively, one may consider a setting where the objective function is evaluated approximately for each sample point, with some high probability of this approximation being accurate, but yet some small probability of a bad approximation. In this case the resulting interpolation/regression model will provide a good approximation with high probability. Note that when computing the function value at the potential new iterate (rather than a sample point) one is still assuming that an accurate value is computed. Relaxing this condition is also a subject for future study.

\paragraph*{Reusing sample points.}

  In a sequential computational setting with expensive function evaluations  it is  efficient to reuse existing sample points in the vicinity of the current iteration. The success of the second method in the example above indicates that sparse models  based on greedy sample sets are useful, even though the sparse recovery properties are unlikely to hold for such sets.
  Hence the random sample models may be dependent in some practical approaches. Investigating the case when the submartingale property holds for such sample sets,
   relaxing the submartingale property in a controlled way, and deriving new convergence results is a subject of our future research.

\subsubsection*{Acknowledgements}

We would like to thank Jose Blanchet and Ramon van Handel for helpful discussions on martingale theory. We also acknowledge Boris Alexeev and Dustin Mixon for interesting discussions on this topic.

\footnotesize

\bibliographystyle{plain}
\bibliography{ref-random_tr}

\begin{thebibliography}{10}

\bibitem{CAudet_JEDennis_2006}
C.~Audet and J.~E. \mbox{Dennis Jr.}
\newblock Mesh adaptive direct search algorithms for constrained optimization.
\newblock {\em SIAM J. Optim.}, 17:188--217, 2006.

\bibitem{UAyaz_HRauhut_2013}
U.~Ayaz and H.~Rauhut.
\newblock Nonuniform sparse recovery with subgaussian matrices.
\newblock {\em ETNA}, 2013, to appear.

\bibitem{ASBandeira_KScheinberg_LNVicente_2012}
A.~S. Bandeira, K.~Scheinberg, and L.~N. Vicente.
\newblock Computation of sparse low degree interpolating polynomials and their
  application to derivative-free optimization.
\newblock {\em Math. Program.}, 134:223--257, 2012.

\bibitem{SCBillups_JLarson_PGraf_2013}
S.~C. Billups, J.~Larson, and P.~Graf.
\newblock Derivative-free optimization of expensive functions with
  computational error using weighted regression.
\newblock {\em SIAM J. Optim.}, 23:27--53, 2013.

\bibitem{RByrd_etal_2012}
R.~Byrd, G.~M. Chin, J.~Nocedal, and Y.~Wu.
\newblock Sample size selection in optimization methods for machine learning.
\newblock {\em Math. Program.}, 134:127--155, 2012.

\bibitem{ECandes_2006}
E.~J. Cand\`es.
\newblock Compressive sampling.
\newblock {\em Proceedings of the International Congress of Mathematicians
  Madrid 2006}, Vol. III, 2006.

\bibitem{ZChen_JJDongarra_2005}
Z.~Chen and J.~J. Dongarra.
\newblock Condition numbers of {G}aussian random matrices.
\newblock {\em SIAM J. Matrix Anal. Appl.}, 27:603--620, 2005.

\bibitem{ARConn_NIMGould_PhLToint_2000}
A.~R. Conn, N.~I.~M. Gould, and \mbox{Ph.} L.~Toint.
\newblock {\em Trust-Region Methods}.
\newblock MPS-SIAM Series on Optimization. SIAM, Philadelphia, 2000.

\bibitem{ARConn_KScheinberg_PhLToint_1997a}
A.~R. Conn, K.~Scheinberg, and \mbox{Ph.} L.~Toint.
\newblock On the convergence of derivative-free methods for unconstrained
  optimization.
\newblock In M.~D. Buhmann and A.~Iserles, editors, {\em Approximation Theory
  and Optimization, Tributes to M. J. D. Powell}, pages 83--108. Cambridge
  University Press, Cambridge, 1997.

\bibitem{ARConn_KScheinberg_PhLToint_1997b}
A.~R. Conn, K.~Scheinberg, and \mbox{Ph.} L.~Toint.
\newblock Recent progress in unconstrained nonlinear optimization without
  derivatives.
\newblock {\em Math. Program.}, 79:397--414, 1997.

\bibitem{ARConn_KScheinberg_PhLToint_1998}
A.~R. Conn, K.~Scheinberg, and \mbox{Ph.} L.~Toint.
\newblock A derivative free optimization algorithm in practice.
\newblock In {\em Proceedings of the 7th AIAA/USAF/NASA/ISSMO Symposium on
  Multidisciplinary Analysis and Optimization, St. Louis, Missouri, September
  2-4}, 1998.

\bibitem{ARConn_KScheinberg_LNVicente_2008a}
A.~R. Conn, K.~Scheinberg, and L.~N. Vicente.
\newblock Geometry of interpolation sets in derivative free optimization.
\newblock {\em Math. Program.}, 111:141--172, 2008.

\bibitem{ARConn_KScheinberg_LNVicente_2008b}
A.~R. Conn, K.~Scheinberg, and L.~N. Vicente.
\newblock Geometry of sample sets in derivative free optimization: {P}olynomial
  regression and underdetermined interpolation.
\newblock {\em IMA J. Numer. Anal.}, 28:721--748, 2008.

\bibitem{ARConn_KScheinberg_LNVicente_2009-paper}
A.~R. Conn, K.~Scheinberg, and L.~N. Vicente.
\newblock Global convergence of general derivative-free trust-region algorithms
  to first and second order critical points.
\newblock {\em SIAM J. Optim.}, 20:387--415, 2009.

\bibitem{ARConn_KScheinberg_LNVicente_2009}
A.~R. Conn, K.~Scheinberg, and L.~N. Vicente.
\newblock {\em Introduction to Derivative-Free Optimization}.
\newblock MPS-SIAM Series on Optimization. SIAM, Philadelphia, 2009.

\bibitem{RDurret_2010}
R.~Durrett.
\newblock {\em Probability: Theory and Examples}.
\newblock Cambridge Series in Statistical and Probabilistic Mathematics.
  Cambridge University Press, Cambridge, fourth edition, 2010.

\bibitem{AEdelman_1988}
A.~Edelman.
\newblock Eigenvalues and condition numbers of random matrices.
\newblock {\em SIAM J. Matrix Anal. Appl.}, 9:543--560, 1988.

\bibitem{GFasano_JLMorales_JNocedal_2009}
G.~Fasano, J.~L. Morales, and J.~Nocedal.
\newblock On the geometry phase in model-based algorithms for derivative-free
  optimization.
\newblock {\em Optim. Methods Softw.}, 24:145--154, 2009.

\bibitem{SGhadimi_GLan_2012}
S.~Ghadimi and G.~Lan.
\newblock Stochastic first- and zeroth-order methods for nonconvex stochastic
  programming.
\newblock Technical report, University of Florida, 2012.

\bibitem{TGKolda_RMLewis_VTorczon_2003}
T.~G. Kolda, R.~M. Lewis, and V.~Torczon.
\newblock Optimization by direct search: {N}ew perspectives on some classical
  and modern methods.
\newblock {\em SIAM Rev.}, 45:385--482, 2003.

\bibitem{JMatyas_1965}
J.~Matyas.
\newblock Random optimization.
\newblock {\em Automation and Remote Control}, 26:246--253, 1965.

\bibitem{JJMore_SMWild_2009}
J.~J. Mor{\'e} and S.~M. Wild.
\newblock Benchmarking derivative-free optimization algorithms.
\newblock {\em SIAM J. Optim.}, 20:172--191, 2009.

\bibitem{YNesterov_2011}
Y.~Nesterov.
\newblock Random gradient-free minimization of convex functions.
\newblock Technical Report 2011/1, CORE, 2011.

\bibitem{YNesterov_2012}
Y.~Nesterov.
\newblock Efficiency of coordinate descent methods on huge-scale optimization
  problems.
\newblock {\em SIAM J. Optim.}, 22:341--362, 2012.

\bibitem{BTPoljak_1987}
B.~T. Polyak.
\newblock {\em Introduction to Optimization}.
\newblock Optimization Software, 1987.

\bibitem{MJDPowell_1994}
M.~J.~D. Powell.
\newblock A direct search optimization method that models the objective and
  constraint functions by linear interpolation.
\newblock In S.~Gomez and J.-P. Hennart, editors, {\em Advances in Optimization
  and Numerical Analysis, Proceedings of the Sixth Workshop on Optimization and
  Numerical Analysis, Oaxaca, Mexico}, volume 275 of {\em Math. Appl.}, pages
  51--67. Kluwer Academic Publishers, Dordrecht, 1994.

\bibitem{MJDPowell_2001}
M.~J.~D. Powell.
\newblock On the {L}agrange functions of quadratic models that are defined by
  interpolation.
\newblock {\em Optim. Methods Softw.}, 16:289--309, 2001.

\bibitem{MJDPowell_2003}
M.~J.~D. Powell.
\newblock On trust region methods for unconstrained minimization without
  derivatives.
\newblock {\em Math. Program.}, 97:605--623, 2003.

\bibitem{MJDPowell_2004}
M.~J.~D. Powell.
\newblock Least {F}robenius norm updating of quadratic models that satisfy
  interpolation conditions.
\newblock {\em Math. Program.}, 100:183--215, 2004.

\bibitem{HRauhut_2010}
H.~Rauhut.
\newblock Compressive sensing and structured random matrices.
\newblock In M.~Fornasier, editor, {\em Theoretical Foundations and Numerical
  Methods for Sparse Recovery}, Radon Series Comp. Appl. Math., pages 1--92.
  2010.

\bibitem{KScheinberg_PhLToint_2010}
K.~Scheinberg and Ph.~L. Toint.
\newblock Self-correcting geometry in model-based algorithms for
  derivative-free unconstrained optimization.
\newblock {\em SIAM J. Optim.}, 20:3512--3532, 2010.

\bibitem{VTorczon_1997}
V.~Torczon.
\newblock On the convergence of pattern search algorithms.
\newblock {\em SIAM J. Optim.}, 7:1--25, 1997.

\bibitem{LNVicente_2013}
L.~N. Vicente.
\newblock Worst case complexity of direct search.
\newblock {\em EURO Journal on Computational Optimization}, 1, 2013.

\bibitem{LNVicente_ALCustodio_2009}
L.~N. Vicente and A.~L. Cust\'odio.
\newblock Analysis of direct searches for discontinuous functions.
\newblock {\em Math. Program.}, 133:299--325, 2012.

\bibitem{SMWild_2008}
S.~M. Wild.
\newblock {MNH: A} derivative-free optimization algorithm using minimal norm
  {H}essians.
\newblock In {\em Tenth Copper Mountain Conference on Iterative Methods}, April
  2008.

\end{thebibliography}

\end{document}